%% file: non_hyperbolic_ergodic_measures.tex
\documentclass[12pt]{amsart}
\usepackage[latin1]{inputenc}
\usepackage{amsmath}
\usepackage{amssymb}
\usepackage{graphicx}
\usepackage{graphics}
\usepackage{amsthm}
\usepackage{amscd}
\usepackage{color}
\usepackage{amsfonts}
\usepackage{fancyhdr}

\newtheorem{theorem}{Theorem}[section]
\newtheorem{mainthm}{Theorem}
\newtheorem*{theorem*}{Theorem}

\newtheorem{proposition}[theorem]{Proposition}
\newtheorem{lemma}[theorem]{Lemma}
\newtheorem{Remark}[theorem]{Remark}
\newtheorem*{definition*}{Definition}

\newtheorem{claim}[theorem]{Claim}
\newtheorem{definition}[theorem]{Definition}
\newtheorem*{Question*}{Question}

\def\N{\mathbb{N}}
\def\Z{\mathbb{Z}}

\def\R{\mathbb{R}}

\def\norm #1{\Vert \,#1\, \Vert\,}

\def\ud{\mathrm{d}}
\def\diff {\operatorname{Diff}}
\def\dim{\operatorname{dim}}

 \def\RR{{\mathbb R}}  
   
 \def\ZZ{{\mathbb Z}}

\def\La{\Lambda}
\def\De{\Delta}

    \def\cT{\mathcal{T}}
\def\cC{\mathcal{C}}   \def\cO{\mathcal{O}} \def\cU{\mathcal{U}}
\def\cD{\mathcal{D}}    \def\cV{\mathcal{V}}
\def\cE{\mathcal{E}}    
\def\cF{\mathcal{F}}   \def\cR{\mathcal{R}}

\keywords{
blender,  robust cycle,
partial hyperbolicity,
ergodic measure,
Lyapunov exponent,
non-hyperbolic measure,
quasi hyperbolic string.
}

\subjclass[2010]{37D30, 37C40, 37C50, 37D25, 37A35}

\begin{document}
\vspace{-2cm}
\title[ non-hyperbolic ergodic measures as the limit of periodic measures]{On the existence of non-hyperbolic ergodic measures as the limit of periodic measures}
\author{Christian BONATTI and Jinhua ZHANG}
\vspace{-2cm}
\maketitle

\begin{abstract}\cite{GIKN} and \cite{BBD1} propose two very different ways for building non hyperbolic measures, \cite{GIKN} building such a measure as the limit of periodic measures
and \cite{BBD1} as the $\omega$-limit set of a single orbit, with a uniformly vanishing Lyapunov exponent. The technique in \cite{GIKN} was essentially used in a generic setting, as the periodic orbits were built by
small perturbations. It is not known if the measures obtained by the  technique in \cite{BBD1} are accumulated by periodic measures.

In this paper we use a shadowing lemma from \cite{G}:
\begin{itemize}\item for getting the periodic orbits in \cite{GIKN} without perturbing the dynamics,
\item for recovering the compact set in \cite{BBD1} with a uniformly vanishing Lyapunov exponent by considering the limit of periodic orbits.
\end{itemize}
As a consequence,
we prove that there exists an open and dense subset $\cU$ of the set of robustly transitive non-hyperbolic diffeomorphisms far from homoclinic tangencies,
such that for any $f\in\cU$, there exists a non-hyperbolic ergodic measure  with full support and approximated by hyperbolic periodic measures.

We also prove that there exists an open and dense subset $\cV$ of the set of  diffeomorphisms exhibiting a robust cycle, such that for any $f\in\cV$,
there exists a non-hyperbolic ergodic measure approximated by hyperbolic periodic measures.
\end{abstract}


\input{intro}

\input{preliminaries}

\input{existence}
\input{proofs}
\input{biblio}

\vskip 2mm

\noindent Christian Bonatti,

\noindent {\small Institut de Math\'ematiques de Bourgogne\\
UMR 5584 du CNRS}

\noindent {\small Universit\'e de Bourgogne, 21004 Dijon, FRANCE}

\noindent {\footnotesize{E-mail : bonatti@u-bourgogne.fr}}

\vskip 2mm
\noindent Jinhua Zhang,

\noindent{\small School of Mathematical Sciences\\
Peking University, Beijing 100871, China}

\noindent {\footnotesize{E-mail :zjh200889@gmail.com }}\\
\noindent and\\
\noindent {\small Institut de Math\'ematiques de Bourgogne\\
UMR 5584 du CNRS}

\noindent {\small Universit\'e de Bourgogne, 21004 Dijon, FRANCE}\\
\noindent {\footnotesize{E-mail : jinhua.zhang@u-bourgogne.fr}}


\end{document}

%% file: intro.tex
\section{Introduction}
In 1960s, R. Abraham and S. Smale \cite{AS} built an example breaking the dream of the density of hyperbolic systems among differentiable systems.
Then people start to search for the robustly non-hyperbolic phenomena. M. Shub \cite{Sh} and R. Ma\~n\'e \cite{M}
built robustly transitive non-hyperbolic examples on $\mathbb{T}^4$ and $\mathbb{T}^3$ respectively. Using ``blender", \cite{BD1}
built more general examples of robustly transitive non-hyperbolic diffeomorphisms on  manifolds of dimension at least three.

 In the celebrated paper \cite{O}, V.Oseledets proves that  any ergodic measure is associated to some real numbers, called \emph{Lyapunov exponents},
 which describe the  asymptotic behavior of the system on the tangent space over a full measure set. An ergodic measure is called \emph{hyperbolic} if all of  its Lyapunov exponents are non-zero.
 Pesin's theory shows that many properties of hyperbolic systems survive for generic point of a hyperbolic ergodic measure.
 Enlightened by this work, people start to study \emph{non-uniformly hyperbolic systems}, that is, systems whose ergodic measures are all hyperbolic.
 There are non-uniformly hyperbolic systems which are not hyperbolic, (see for instance
\cite{CLR} which build a non-uniformly hyperbolic surface diffeomorphism  exhibiting homoclinic tangency). However, it has been recently announced that  such
non-hyperbolic but non-uniformly hyperbolic systems cannot be $C^1$-robust (see \cite{CCGWY}).

Comparing with hyperbolic situation,  it was natural to ask whether the non-uniform hyperbolic systems are dense among the differentiable systems.
The answer is ``No'': \cite{KN} builds a $C^1$-open set of diffeomorphisms $f$ having  a non-hyperbolic ergodic measure.
The example in \cite{KN} uses a precise global setting (partially hyperbolic system obtained as small perturbations of skew product).
Recently,  \cite{BBD1} shows that the $C^1$-robust existence of non-hyperbolic ergodic measure is open and dense among diffeomorphisms having a \emph{robust cycle}
(that is, hyperbolic sets of different indices whose stable manifolds of both of them intersect robustly the unstable manifold of the other).

 The approaches of \cite{KN} and \cite{BBD1} are very different:
 \begin{itemize}
 \item \cite{KN} builds the non-hyperbolic measure as being the weak$*$-limit of a sequence $\{\mu_n\}_{n\in\N}$  of \emph{periodic measures } i.e. $\mu_n$  is the measure supported on a periodic orbit $\cO_n$.
 The orbits $\cO_n$ follow \emph{a criterion from \cite{GIKN}}, explicitly stated in \cite[Lemma 2.5]{BDG},  ensuring that the limit of the $\{\mu_n\}_{n\in\N}$ is ergodic.
 \item \cite{BBD1} builds a partially hyperbolic compact invariant set of non-vanishing topological entropy so that every point has
 well defined and vanishing center Lyapunov exponent. This compact set is built as being the $\omega$-limit set of a single point which satisfies a criterion called
 \emph{controlled at any scale}.
 \end{itemize}
 It is not understood up to now if the measures built following the criterion in \cite{GIKN} may have positive entropy,
 and in contrast it is not clear if the non-hyperbolic ergodic measures supported on the compact set built in \cite{BBD1} are accumulated by periodic measures.

 In this paper, we analyze the periodic orbits in a neighborhood of the robust cycle defined in \cite{BBD1},
 and which are less and less hyperbolic: some of their Lyapunov exponents tends to $0$. \cite{DG,BDG} already considered such sequence of periodic orbits in a neighborhood of a robust cycle,
 but they  build the periodic orbits by performing perturbations of the dynamics, so that their conclusion was in term of $C^1$-generic diffeomorphisms.
 Here we do not perturb the dynamics: under some open and dense geometric setting, we prove the existence of such periodic orbits by using the shadowing lemma in \cite{G}.
 Let us present roughly our results:

 \begin{itemize}
 \item On one hand, combining the shadowing lemma by \cite{G} and the \cite{GIKN} criterion
 we build  a sequence of hyperbolic periodic measures whose weak$*$-limit is a non-hyperbolic ergodic measure (Theorem~\ref{thmD}).
 As a consequence, Theorems~\ref{thm.non-hyperbolic ergodic measure} and~\ref{thm.full support} prove that, for an open and dense subset of robustly
 transitive partially hyperbolic  but non-hyperbolic systems far from homoclinic tangencies, there exists a non-hyperbolic ergodic measure which
 has full support and  is the weak$*$-limit of a sequence of hyperbolic periodic measures satisfying the \cite{GIKN} criterion.
  \item On the other hand, in Theorem~\ref{thmC}, using the \emph{controlled at any scale criterion} and the shadowing lemma by \cite{G},
  we recover the compact invariant set $K'_f$ in \cite{BBD1}, with well defined and vanishing center Lyapunov exponent,
and we build a sequence of hyperbolic  periodic orbits $\cO_n$ whose center Lyapunov exponents tend  to zero, and so that:
\begin{itemize}
\item every weak$*$-limit measure $\mu$  of the  measure $\mu_n$ supported on $\cO_n$   is supported on $K'_f$;
\item the Hausdorff limit $K_f$ of the sequence $\{\cO_n\}_{n\in\N}$ contains $K'_f$.
\item any ergodic measure supported on $K_f$ is either supported on $K'_f$ or is a (unique) periodic measure.
\end{itemize}
 \end{itemize}

We are now ready to present our precise setting and the statements of our results. We start with the consequences of our results in a global setting.
Then we will present the more technical results in the local setting, which are the heart of our work.

\subsection{Results  in the global setting}

Let $\cT(M)$ be the subset of $\diff^1(M)$ such that for any $f\in\cT(M)$, we have that
\begin{itemize}
\item  $f $ is robustly transitive;
\item $f$ admits a partially hyperbolic splitting of the form
$$TM=E^s\oplus E_1^c\oplus\cdots\oplus E^c_k\oplus E^{u}$$
satisfying that $\dim(E_1^c)=\cdots=\dim(E^c_k)=1$.
\item there exist two hyperbolic periodic orbits of indices $\dim(E^s)$ and $\dim(E^s)+k$ respectively.
\end{itemize}

We denote by $\dim(E^s)=i_0$ and $\dim(M)=d$.
\begin{Remark}
\begin{itemize}\item By definition, $\cT(M)$ is an open subset of $\diff^1(M)$.
\item  By robust transitivity and  \cite[Theorem 1]{ABCDW}, there exists an open and dense subset $\cT_p(M)$ of $\cT(M)$ such that for any $f\in\cT_p(M)$ and any $j=0,\cdots, k$,  there exists a  hyperbolic periodic orbit of index $i_0+j$.
\end{itemize}
\end{Remark}

For any $f\in\cT(M)$ and any $f$-ergodic measure $\mu$, we denote by $\lambda^c_i(\mu)$ the Lyapunov exponent of $\mu$ along the bundle $E^c_i$, for any $i=1,\cdots, k$.
As the bundle $E^c_i$ is one dimensional,  one has  that
$$\lambda^c_i(\mu)=\int\log\norm{Df|_{_{E^c_i}}}\ud\mu.$$

\begin{theorem}\label{thm.non-hyperbolic ergodic measure} There exists an open and dense subset $\tilde{\cT}(M)$ of $\cT(M)$, such that for any $f\in\tilde{\cT}(M)$, there exist $k$ non-hyperbolic ergodic measures $\mu_1,\cdots,\mu_k$ such that for any $i=1,\cdots,k$, we have that
\begin{itemize}
\item $\lambda^c_i(\mu_i)=0;$
\item The support of each $\mu_i$ is the whole manifold $M$;
 \item $\mu_i$ is the weak$*$-limit of a sequence of hyperbolic periodic measures of index $i_0+i$;
\end{itemize}
\end{theorem}
Note that in the case $k=1$,  the existence of non-hyperbolic ergodic measure with full support  is also announced in \cite{BBD2}, but the proof is quite different.

We denote by $\cU(M)$ the set of robustly transitive $C^1$ diffeomorphisms far from homoclinic tangencies. Hence, $\cU(M)$ is an open set of $\diff^1(M)$.

It is clear, by   definition, that $\cT(M)$ is a subset of $\cU(M)$. Indeed, by Theorem D in \cite{BDPR}, the set $\cT(M)$
is an open and dense subset of $\cU(M)$. As a consequence of Theorem \ref{thm.non-hyperbolic ergodic measure}, one gets the following straightforward result:
\begin{mainthm}\label{thm.full support}
There exists an open and dense subset $\cV(M)$ of $\cU(M)$,
 such that for any $f\in\cV(M)$, there exists a non-hyperbolic ergodic measure as the weak$*$-limit of a sequence of periodic measures,  whose support is the whole manifold.
\end{mainthm}

 \subsection{Results in the semi-local setting of  robust cycles}
In this paper, we consider diffeomorphisms having a robust cycle.
 The precise setting is defined in \cite{BBD1} and  is called
\emph{flip flop configuration},  whose definition uses many other notions defined in   Section~\ref{s.Preliminaries}.  For this reason, we state our results without
too technical definitions.

Consider the open set $\cR\cC$(M)  of diffeomorphism  $f$ presenting  a robust cycle between a transitive hyperbolic set $\La_f$ of $s$-index $i$ and
hyperbolic periodic point $q_f$ of $s$-index $i+1$.
\cite{BBD1} built an open  and dense subset $\widetilde{\cR\cC}(M)$ (that is, the set of diffeomorphisms with a \emph{split flip-flop configuration}) in $\cR\cC(M)$ so that every $f\in \widetilde{\cR\cC}(M)$ admits a point $x_f$ whose $\omega$-limit set $ \omega(x_f)$ has the following properties:
\begin{itemize}
 \item the topological entropy of $ \omega(x_f)$ is positive;
 \item $ \omega(x_f)$ is partially hyperbolic with $1$-dimensional center bundle;
 \item every $y\in  \omega(x_f)$ has  well defined and vanishing center Lyapunov exponent.
\end{itemize}
This result contrasts with the procedure in \cite{GIKN,KN}  which build non-hyperbolic measure as the limit of periodic orbits, in a specific global setting
(skew product of a hyperbolic dynamics by diffeomorphisms of  the circle). In particular, it is not clear a priori if the (non-hyperbolic) measures supported on the compact set $ \omega(x_f)$ built in
\cite{BBD1} are accumulated by periodic orbits.

Our first result consists in showing that the controlled at any scale criterion can be used with periodic orbits:
the orbits follow this controlled at any scale criterion out of a very small orbit segment, whose weight in the corresponding periodic
measure tends to $0$; this small orbit segment is used for closing the orbit by using a shadowing lemma in \cite{G}.

\begin{mainthm} \label{thmC} With the notations above,  for any $f\in\widetilde{\cR\cC}(M)$, there exists a  sequence of hyperbolic periodic orbits $\{\gamma_n\}$
homoclinically related to the orbit of $q_f$,  whose center Lyapunov exponent tends to zero, and which converges for the Hausdorff distance to a compact invariant set
$K_f$ such that:
  \begin{itemize}\item $ q_f\in  K_f$;
  \item the set  $K_f$ is partially hyperbolic with $1$-dimensional center bundle;
  \item there exists a non-empty compact invariant set $K^{\prime}_f\subset K_f$ such that any point in $K^{\prime}_f$
  has well defined and vanishing center Lyapunov exponent.
      \item  the topological entropy of $K^{\prime}_f$ is positive: $h_{top}{f|_{K^{\prime}_f}}>0$;
      \item For any $x\in K_f\backslash(K^{\prime}_f\cup\mathcal{O}_{q_f})$,
      we have either $\omega(x)\subset K ^{\prime}_f$ and $x\in W^{u}(\mathcal{O}_{q_f})$ or $\alpha(x)\subset K^{\prime}_f$  and  $x\in W^{s}(\mathcal{O}_{q_f})$.
  \end{itemize}
  \end{mainthm}

  \begin{Remark}
  \begin{enumerate}

  \item The invariant compact set $K^{\prime}_f$  built here is indeed    the compact set built in \cite{BBD1}
  (see Remark~\ref{r.choice of Kf} in Section~\ref{s.proof of theorem B}), that is,  the $\omega$-limit set of a
point $x\in W^u(q_f)$ whose positive orbit is controlled at any scale.
   \item Any ergodic measure supported on $K_f$ is either the Dirac measure on $\mathcal{O}_{q_f}$ or a non-hyperbolic ergodic measure.

\item Any limit measure of the periodic measures $\delta_{\gamma_n}$ (supported on $\gamma_n)$ is supported on $K^{\prime}_f$:
  this follows from the fact that the weight given by $\delta_{\gamma_n}$ to
  a neighborhood of  $\cO_{q_f}$ tends to $0$.

 \item As the periodic orbits in Theorem~\ref{thmC} are all homoclinically related, one gets that, for  any limit measure $\mu$  of $\{\delta_{\gamma_n}\} $,
  the whole probability segment $\{t\mu+(1-t)\delta_{\cO_{q_f}}, t\in[0,1]\}$ is accumulated by periodic measures, where $\delta_{\cO_{q_f}}$ denotes
  the periodic measure supported on   $\cO_{q_f}$.
   \end{enumerate}
  \end{Remark}

  Our second local result consists in showing that the criterion from \cite{GIKN} applies for any $f\in \widetilde{\cR\cC}(M)$,
  ie. the diffeomorphisms with a split-flip-flop configuration.  As a consequence,  we get:

  \begin{mainthm}\label{thmD} For any $f\in \widetilde{\cR\cC}(M)$, there is a partially hyperbolic set $\tilde \La_f$
  (with $1$-dimensional center bundle) and a sequence of periodic orbits $\{\mathcal{O}_{n}\}\subset \tilde{\Lambda}_f$ such that:
  \begin{itemize}\item  the center Lyapunov exponent of $\mathcal{O}_{n}$ tends to zero;
  \item the orbits $\cO_{n}$ satisfy the \cite{GIKN} criterion. As a consequence, one has that  the Dirac measure $\delta_{\mathcal{O}_{n}}$ converges to a non-hyperbolic ergodic measure whose support is the Hausdorff limit of the orbits $\{\cO_n\}$.
  \end{itemize}
  \end{mainthm}
  \vskip 5mm
{\bf Acknowledgment: }
We would like to  thank Lorenzo D\'iaz for useful comments.

Jinhua Zhang  would like to thank Institut de Math\'ematiques de
Bourgogne for hospitality and China Scholarship Council (CSC) for financial support (201406010010).

%% file: preliminaries.tex
\section{Preliminaries}\label{s.Preliminaries}
In the whole paper, we assume that $M$ is compact Riemannian manifold.

  In this section, we will collect some notations and some results that we need.  We start by recalling very classical notions, as hyperbolic basic set and  dominated splitting.
  Then we recall our main (more recent) tools.  More precisely, our results consist in applying  four tools in a very specific setting.  The tools are:
  \begin{itemize}
  \item a criterion by \cite{GIKN} for ensuring that a limit measure of measures supported on periodic orbits is ergodic.

  \item a shadowing lemma due to S. Liao \cite{ Liao1} and S. Gan in \cite{G} : this will allow us to prove the existence of  periodic orbits with a prescribed itinerary.

   \item a criterion in \cite{BBD1} (called \emph{controlled  at any scale}) for controlling averages along an orbit. We will apply it here in the partially hyperbolic setting
   for getting a vanishing center Lyapunov exponent.
   \item an abstract dynamical system called \emph{flip flop family} which will be our machinery for producing the orbits on which we can apply the three tools above.
  \end{itemize}
  Our setting will be a specific   robust cycle defined in \cite{BBD1} and called \emph{flip flop configuration}. The main interest of the flip flop configuration is that they appear
  open and densely in the setting of robust cycle, and they provide flip flop families.
\subsection{Dominated splitting,  partial hyperbolicity and hyperbolicity}
 Let us recall that a $Df$-invariant splitting $T_{K}M=E\oplus F$ on a compact $f$-invariant set $K$ is a \emph{dominated splitting}, if there exist
 $\lambda\in(0,1)$ and a metric $\|\cdot\|$ such that for any $x\in K$, we have
 \begin{displaymath} \norm{Df |_{E(x)}}\norm{Df^{-1}|_{F(f(x))}}^{-1}< \lambda.
  \end{displaymath}

  We say that a $Df$-invariant splitting $T_{K}M=E^{s}\oplus E^{c}\oplus E^{u}$ on a compact $f$-invariant set $K$ is a \emph{partially hyperbolic splitting},
  if $E^{s}$ and $E^{u}$ are  uniformly contracting and expanding respectively, and  the splittings $E^{s}\oplus (E^{c}\oplus E^{u})$ and $(E^{s}\oplus E^{c})\oplus E^{u}$ are dominated splittings.

A much stronger splitting is called \emph{hyperbolic splitting}. Recall that a $Df$-invariant splitting $T_{K}M=E^s\oplus E^{u}$ on a compact invariant set
$K$ is a \emph{hyperbolic splitting}, if  $E^s$ is uniformly contracting and $E^{u}$ is uniformly expanding under $Df$, and $K$ is called \emph{hyperbolic set}.
A hyperbolic set $K$ is called a \emph{hyperbolic basic set} if $K$ is transitive and there exists an open neighborhood $U$ of $K$ such that $K$ is the \emph{maximal invariant set in $U$},  that is
$$K=\cap_{i\in\Z}f^{i}(U).$$
\subsection{Center Lyapunov exponent of ergodic measures supported on a partially hyperbolic set}
Let $K$ be an $f$ invariant compact set admitting a partially hyperbolic splitting of the form $T_{K}M=E^s\oplus E^c\oplus E^{u}$ such that $\dim(E^c)=1$. We denote by $\dim(E^s)=i$.

For any ergodic measure $\mu$ supported on $K$, the center Lyapunov exponent of $\mu$ is defined as:
$$\lambda^c(\mu)=\int\log\norm{Df|_{E^c}}\ud\mu.$$

Let $\mu$ be an ergodic measure supported on $K$. $\mu$ is called a non-hyperbolic ergodic measure if we have $\lambda^c(\mu)=0$, and
$\mu$ is called a hyperbolic ergodic measure of index $i$ (resp. $i+1$) if we have that $\lambda^c(\mu)>0$ (resp. $\lambda^c(\mu)<0$).

\subsection{Homoclinic class}
\begin{definition} Given $f\in\diff^1(M)$. Let  $\mathcal{O}_p$ and $\mathcal{O}_q$ be two hyperbolic periodic orbits of $f$. We say that  $\mathcal{O}_p$ and $\mathcal{O}_q$ are \emph{homoclinically related}, if there exist two transverse intersections  $x\in W^{s}(\cO_p)\cap W^{u}(\cO_q)$ and $y\in W^{u}(\cO_p)\cap W^{s}(\cO_q)$.
\end{definition}
We denote by $P(f)$ the set of hyperbolic periodic orbits of the diffeomorphism $f$. Let $\mathcal{O}_p\in P(f)$,  the homoclinic class of $\mathcal{O}_p$ is defined as:
$$H(p,f):=\overline{\{\mathcal{O}_q\in P(f)|\textrm{ $\mathcal{O}_q$ is  homoclinically related to $\mathcal{O}_p$}\}}.$$

Given two hyperbolic periodic orbits $\mathcal{O}_p$ and $\mathcal{O}_q$. Let $U$ be an open  neighborhood of $\mathcal{O}_p\cup\mathcal{O}_q$. We say that $\mathcal{O}_p$ and $\mathcal{O}_q$ are \emph{homoclinically related inside $U$}, if there exist two transverse intersections  $x\in W^{s}(\cO_p)\cap W^{u}(\cO_q)$ and $y\in W^{u}(\cO_p)\cap W^{s}(\cO_q)$ such that $\overline{Orb(x)}\cup \overline{Orb(y)}\subset U.$

  \subsection{A criterion for ergodicity of the limit measure of periodic measures}
  In this subsection, we  state a criterion ensuring that a sequence of periodic measures converges to an ergodic measure.
  This criterion is firstly used in \cite{GIKN,KN} and then developed in \cite{BDG} for building  non-hyperbolic ergodic measures as limit of periodic measures whose center Lyapunov exponents tend to $0$.
  \begin{definition}\label{good for}
  Let $(X, \ud)$ be a compact metric space and $f:X\mapsto X$ be a continuous map.
  Fix $\epsilon>0$ and $\kappa\in(0,1)$.  Let $\gamma_1$ and $\gamma_2$ be two periodic orbits of $f$.
   Then, the periodic orbit $\gamma_{1}$ is said to be \emph{$(\epsilon,\kappa)$ good }for the  periodic orbit $\gamma_{2}$ if the followings hold:
  \begin{itemize}
  \item there exist a subset $\gamma_{1,\epsilon}$ of $\gamma_{1}$ and a projection $\xi :\gamma_{1,\epsilon}\rightarrow\gamma_{2}$ such that $$\ud(f^{i}(y),f^{i}(\xi(y)))<\epsilon,$$ for every $y\in\gamma_{1,\epsilon}$ and every $j=0,1,\ldots,\pi(\gamma_{2})-1$;
      \item the proportion of $\gamma_{1,\epsilon}$ in $\gamma_1$ is larger that $\kappa$.  In formula:
      $$\frac{\#\gamma_{1,\epsilon}}{\pi(\gamma_{1})}\geq\kappa.$$
      \item  the cardinal of the preimage $\#\xi^{-1}(x)$  is the same for all $x\in\gamma_{2}$.
  \end{itemize}
  \end{definition}
  We can now state the \emph{\cite{GIKN} criterion}:
  \begin{lemma}\label{limit} \cite[Lemma 2.5]{BDG}
   Let $(X, \ud)$ be a compact metric space and $f:X\mapsto X$ be a homeomorphism. Let $\{X_{n}\}$ be a sequence of periodic orbits
  whose periods tend to infinity. Let $\mu_{n}$ denote  the Dirac measure of $X_{n}$.

  Assume that the orbit $X_{n+1}$ is $(\epsilon_{n},  \kappa_{n})$ good for $X_{n}$, where   $\epsilon_{n}>0$ and
   $0<\kappa_{n}<1$ satisfy
  \begin{displaymath} \sum_{n}\epsilon_{n}<\infty \textrm{ and } \prod_{n}\kappa_{n}>0.
  \end{displaymath}

  Then the sequence $\{\mu_{n}\}$ converges to an ergodic measure $\nu$ whose support is given by
  \begin{displaymath} supp\,\nu=\cap_{n=1}^{\infty}\overline{\cup_{k=n}^{\infty}X_{k}}.
  \end{displaymath}
\end{lemma}

\subsection{A shadowing lemma}
In this paper we don't construct any periodic orbits by perturbing the dynamics; we just find out these periodic orbits.
The way we use to detect these periodic orbits is a shadowing lemma which is firstly given by S.  Liao \cite{Liao1} and is developed by S. Gan\cite{G}.

    Let $\Lambda$ be an $f$-invariant compact set.
    Assume that there exists a $Df$-invariant continuous splitting $T_{\Lambda}M=E\oplus F$.
      For any $\lambda<0$, an orbit segment $\{x,n\}:=\{x,\ldots,f^{n}(x)\}$ contained in $\Lambda$ is called a \emph{$\lambda$-quasi hyperbolic string}, if the followings are satisfied:
      \begin{itemize}
       \item {\bf Uniform contraction of $E$ by $Df$, from $x$ to $f^n(x)$:}
       $$\frac{1}{k}\sum_{i=0}^{k-1}\log{\|Df|_{E(f^{i}(x))}\|}\leq\lambda,$$ for every  $k=1,\cdots,n$;
       \item {\bf Uniform contraction of $F$ by $Df^{-1}$, from $f^n(x)$ to $x$}
       $$\frac{1}{n-k}\prod_{i=k}^{n-1}\log{m(Df|_{F(f^{i}(x))})}\geq -\lambda,$$ for every $k=0,\cdots,n-1$;
       \item {\bf Domination of $E$ by $F$ from $x$ to $f^n(x)$ }
       $$\log{\|Df|_{E(f^{i}(x))}\|}-\log{m(Df|_{F(f^{i}(x))})}\leq 2\lambda,$$  for every $i=0,\cdots,n-1$.
      \end{itemize}

 \begin{Remark}
From the definition, we can easily check  that a $\lambda$-quasi hyperbolic string is also a $\frac{\lambda}{2}$-quasi hyperbolic string.
\end{Remark}

\begin{definition}
 Consider   $d>0$ and $\lambda<0$. Let $\{x_i\}_{i\in\Z}$ be a sequence of points in $\Lambda$ and $\{n_i\}_{i\in\Z}$
 be a sequence of positive integers.
 We say that the sequence of orbit segments $\big\{\{x_i,n_i\}\big\}_{i\in\Z}$ is \emph{a $\lambda$-quasi hyperbolic $d$-pseudo orbit}
 if for any $i$, we have:
 \begin{itemize}\item  $\ud(f^{n_i}(x_i),x_{i+1})\leq d$,
 \item  the orbit segment $\{x_i,n_i\}$ is a $\lambda$-quasi hyperbolic string.
  \end{itemize}

 \end{definition}
 We say that a $\lambda$-quasi hyperbolic $d$-pseudo orbit  $\big\{\{x_{i},n_{i}\}\big\}_{i\in\Z}$  is \emph{periodic},
 if  there exists a positive integer $m$ such that $n_{i+m}=n_{i}$ and $x_{i+m}=x_{i}$ for all $i$.
 Then, assuming that $m$ is the smallest such positive integer, the sum $\sum_{i=0}^{m-1} n_i$ is  \emph{the period} of the pseudo orbit.

 \begin{definition} Let $\{x_i\}_{i\in\ZZ}$ be a sequence of points and $\{n_i\}_{i\in\ZZ}$ be a sequence of strictly positive integers.
 We define
 \begin{displaymath}T_i=\left\{\begin{array}{ll}
 0 & \textrm{if $i=0$}\\
 n_0+\cdots+n_{i-1} & \textrm{if $i>0$}\\
 -n_{-i}-\cdots-n_{-1} & \textrm{if $i<0$}
 \end{array}\right.
 \end{displaymath}

Let $\epsilon>0$, we say that the orbit of a point $x$ \emph{$\epsilon$-shadows $\big\{\{x_i,n_i\}\big\}_{i\in\Z}$}  if for any
$i\in\ZZ$ and $T_i\leq j\leq T_{i+1}-1$, we have that $$\ud(f^{j}(x),f^{j-T_i}(x_i))<\epsilon.$$
 \end{definition}

 \begin{lemma} \label{shadow}\cite{G} $[$Shadowing lemma for quasi hyperbolic pseudo orbit$]$  Assume that $\Lambda$ is an $f$-invariant compact set
  and there exists an $f$-invariant continuous splitting $T_{\Lambda}M=E\oplus F$.

  Then, for any $\lambda<0$, there exist $L>0$ and $d_{0}>0$ such that for any $d\in(0,d_{0}]$ and any
  $\lambda$-quasi hyperbolic $d$-pseudo orbit  $\big\{\{x_{i},n_{i}\}\big\}_{i\in\Z}$,
  there exists a point $x$ whose orbit $L\cdot d$ shadows $\big\{\{x_{i},n_{i}\}\big\}_{i\in\Z}$ .

  If moreover the quasi hyperbolic pseudo-orbit $\big\{\{x_{i},n_{i}\}\big\}_{i\in\Z}$  is periodic,
  then the point $x$ can be chosen to be periodic with the same period.
  \end{lemma}
\subsection{Plaque family,  hyperbolic time  and estimate on the  size of the invariant manifold}
Let $T_{\Lambda}M=E\oplus F$ be a dominated splitting over a compact $f$-invariant set $\Lambda$.
We denote by $\dim(E)=i$ and $\dim(M)=n$.
  Let  $\mathbb{D}^{i}$ be the $i$-dimensional unit disc and $\mathbb{D}^{n-i}$ be the $(n-i)$-dimensional unit disc.
  In addition, we denote by $Emb^{1}(\mathbb{D}^{i},M)$ the space of $C^1$-embedding maps from $\mathbb{D}^{i}$ to $M$ and  $Emb^{1}(\mathbb{D}^{n-i},M)$ the space of $C^1$-embedding maps from $\mathbb{D}^{n-i}$ to $M$.

 Like the situation of hyperbolic set, for a compact set with dominated splitting, there also exist invariant manifolds   due to \cite{HPS}. To be precise:  \begin{lemma}\label{l.plaque family}
 Under the assumption and notation above, there exist two families of continuous maps
 $$\Phi^{cs}: \Lambda\rightarrow Emb^{1}(\mathbb{D}^{i},M) \textrm{ and }\Phi^{cu}: \Lambda\rightarrow Emb^{1}(\mathbb{D}^{n-i},M).$$
   We denote  by  $\mathcal{W}^{cs}(x)=\Phi^{cs}(x)(\mathbb{D}^{i})$  and $\mathcal{W}^{cu}(x)=\Phi^{cu}(x)(\mathbb{D}^{n-i})$, then  the following properties hold:
  \begin{itemize}\item $T_{x}\mathcal{W}^{cs}(x)=E_{x}$ and $T_{x}\mathcal{W}^{cu}(x)=F_{x}$;
  \item For any $\delta_{1}\in(0,1)$, there exists $\delta_{2}>0$ such that:
$$f({W}_{\delta_{2}}^{cs}(x))\subset {W}_{\delta_{1}}^{cs}(f(x)) \textrm{ and }
f^{-1}({W}_{\delta_{2}}^{cu}(x))\subset {W}_{\delta_{1}}^{cu}(f^{-1}(x)), \textrm{ for any $x\in\Lambda$}.$$

 \end{itemize}
 We call $\{W^{cs}(x)\}_{x\in\Lambda}$ and $\{W^{cu}(x)\}_{x\in\Lambda}$
  the \emph{plaque families} of  the dominated splitting $E\oplus F$.
 \end{lemma}

\begin{definition} Let $\La$ be a compact invariant set admitting  a dominated splitting $T_{\La}M=E\oplus F$. Given $\lambda<0$.
A point $x\in\Lambda$ is called a \emph{$(\lambda, E)$ hyperbolic  time} if we have the following:
$$\sum_{i=0}^{n-1}\log\norm{Df|_{E(f^{i}(x))}}\leq n\cdot \lambda, \textrm{ for any  integer $n\geq 1$.}$$
Similarly, we can define the \emph{ $(-\lambda ,F)$-hyperbolic time} which is a $(\lambda, F)$ hyperbolic time for $f^{-1}$.
\end{definition}
By Lemma \ref{l.plaque family}, we can  fix a plaque family $W^{cs}$ corresponding to the bundle $E$.
The following lemma guarantees the existence of stable manifold at the $(\lambda, E)$ hyperbolic time. The proof is classical (see for instance \cite[Section8.2]{ABC}).
\begin{lemma}\label{lp} For any $\lambda<0$, there exists $\eta>0$  such that for any $(\lambda, E)$ hyperbolic time $x$, we have that  the disc $W^{cs}_{\eta}(x)$ is contained in the stable manifold of $x$.
  \end{lemma}
  \begin{Remark} Similar result holds for  $(-\lambda, F)$ hyperbolic time.
  \end{Remark}
  To find the  $(\lambda, E)$ hyperbolic time, we need the following well known result:
  \begin{lemma}\label{l.Pliss lemma}\cite{P} Given a number $A$. Consider a sequence of numbers  $a_{1},\cdots,a_{n}$  bounded from above by $A$. Assume that there exists a number $c<A$ such that $\sum_{i=1}^{n}a_{i}\geq n\cdot c$.

    Then for any number $c^{\prime}<c$, there exist $l$ integers $t_{1},\cdots,t_{l}\subset [1,n]$ satisfying  that:
    $$\sum_{i=j}^{t_{k}}a_{i}\geq (t_{k}-j+1)c^{\prime},\textrm{ for any $k=1,\cdots,l$ and any $j=1,\cdots,t_{k}$}.$$
    Moreover, we have  $\frac{l}{n}\geq\frac{c-c^{\prime}}{A-c^{\prime}}$.
  \end{lemma}

  Let $p$ be a periodic point and $\lambda$ be a negative number. Assume that there exists a $Df$ invariant splitting $T_{\cO_p}M=E\oplus F$ over the orbit $\cO_p$. The point $p$ is called \emph{a $\lambda$ bi-hyperbolic time}  if   for any $k=1,\cdots,\pi(p)$, we have that
 $$\frac 1k \sum_{i=0}^{k-1}\log\norm{Df|_{E(f^{i}(p))}}\leq\lambda$$
and
 $$\frac 1k \sum_{i=0}^{k-1}\log\norm{Df^{-1}|_{F(f^{-i}(p))}}\leq\lambda.$$

 The following classical lemma gives the existence of bi-hyperbolic time ( see for instance \cite[Lemma 2.21]{Wa}).
 \begin{lemma}\label{l.bi-hyperbolic time}Let $f\in\diff^1(M)$ and  $p$ be a periodic point. Assume that there exists a $Df$ invariant splitting $T_{\cO_p}M=E\oplus F$ and a number $\lambda<0$ satisfying that
 \begin{itemize} \item  $$\frac {1}{\pi(p)} \sum_{i=0}^{\pi(p)-1}\log\norm{Df|_{E(f^{i}(p))}}\leq\lambda$$
\item
 $$\frac{1}{\pi(p)} \sum_{i=0}^{\pi(p)-1}\log\norm{Df^{-1}|_{F(f^{-i}(p))}}\leq\lambda.$$
 \item $\log\norm{Df|_{E(f^{i}(p))}}+\log\norm{Df^{-1}|_{F(f^{i+1}(p))}}\leq 2\cdot\lambda, \textrm{ for any integer $i$.}$
 \end{itemize}

 Then for any $\lambda^{\prime}\in(\lambda,0)$, there exists a point $q\in \cO_p$ such that $q$ is a $\lambda^{\prime}$ bi-hyperbolic time.
 \end{lemma}
 \subsection{Control of averages at any scale}

    In this subsection, we restate a criterion given in \cite{BBD1} of the existence of zero average   for a continuous function along an orbit.
 In this section,  let   $(X,\ud)$ be  a metric space, $K\subset X$ be a compact subset, $f:X\mapsto X$ be  a homeomorphism  and $\varphi:K\rightarrow \R$ be a continuous function.

    \begin{definition} Given $\beta>0$, $t\in\N$  and $T\in\N^{+}\cup\{+\infty\}$, we say that a point $x\in K$ is \emph{ $(\beta,t,T)$-controlled},
    if $f^{i}(x)\in K$ for $0\leq i< T$ and there exists a subset $\mathcal{P}\subset\N$  such that
    \begin{itemize}\item $0\in\mathcal{P}$ and $T = sup (\mathcal{P})$;
    \item if $k<l$ are two consecutive elements in $\mathcal{P}$,  then  we have
    $$l-k\leq t \textrm{  and  }\frac{1}{l-k}\big|\sum_{i=0}^{l-k-1}\varphi\big(f^{i+k}(x)\big)\big|\leq \beta.$$

    \end{itemize}
    \end{definition}

    A point $x\in K$ is \emph{controlled at any scale} if there exist   monotone sequences $(t_{i})_{i}$ of natural numbers and $(\beta_{i})_{i}$
    of positive numbers, with $t_{i}\rightarrow +\infty$ and $\beta_{i}\rightarrow 0$, such that $x$ is $(\beta_{i},t_{i}, +\infty)$-controlled for every $i$.
    Note that this implies that the $\omega$-limit set $\omega(x)$ is contained in $K$.

    Denote by $\varphi_{n}(x):=\sum_{i=0}^{n-1}\varphi(f^{i}(x))$, for $x\in\cap_{i=0}^{n-1}f^{-i}(K)$. In particular, if $x$ is controlled at any scale,
    its positive orbit remains in $K$ so that $\varphi_n$ is defined and continuous on the closure of this positive orbit.
    Now, for the points controlled at any scales, we have the following property:
    \begin{lemma}\label{l.control at any scale}\cite[Lemma 2.2]{BBD1} If $x\in K$ is controlled at any scale, then every point $y\in\omega(x)$ satisfies
      $$\lim_{n\rightarrow\infty}\frac{1}{n}\varphi_{n}(y)=0.$$
      Moreover, the limit is uniform over the $\omega$-limit set $\omega(x)$.
    \end{lemma}
\subsection{Flip flop family and the control at any scale}
Let $(X,\ud)$ be a compact metric space and $f:X\mapsto X$ be a continuous map. Let $K$ be a compact subset of $X$ and $\varphi:K\mapsto\R$ be a continuous function.
    \begin{definition}\label{def.flip-flop family}$(Flip-flop\, family)$
    A \emph{flip-flop family}, associated to the continuous function $\varphi$, is a family $\mathfrak{F}$ of compact subsets of $K$ with uniformly bounded diameters
    that splits as $\mathfrak{F}=\mathfrak{F}^{+}\cup\mathfrak{F}^-$ into two disjoint families satisfying:
    \begin{enumerate}\item There exists a constant $\alpha$ such that for any $D^+\in\mathfrak{F}^+,\,D^-\in\mathfrak{F}^-$ and any points
    $x\in D^+, y\in D^-$, we have $\varphi(x)>\alpha>0>-\alpha>\varphi(y)$;
    \item For any $D\in\mathfrak{F}$, there exist two subsets $D_+,\, D_-$ of $D$ such that $f(D_+)\in\mathfrak{F}^+$ and $f(D_-)\in\mathfrak{F}^-$;
        \item There exists a constant $\lambda>1$ such that if $x,y$ belong to the same element $D_0$ of $\mathfrak{F}$ and if $f(x)$ and $f(y)$ belong also to
        the same element $D_1$ of $\mathfrak{F}$ then
      $$\ud(f(x),f(y))\geq\lambda \,\ud(x,y).$$
    \end{enumerate}
    \end{definition}
    An important property of flip-flop family is the following:
    \begin{lemma}\cite[Theorem 2.1]{BBD1} For any $D\in\mathfrak{F}$, there exists a point $x\in D$ such that $x$ is controlled at any scale with respect to $\varphi$.
    \end{lemma}

    By Definition \ref{def.flip-flop family}, we can iterate any element of $\mathfrak{F}$ and its image contains an element of $\mathfrak{F}$.
    This leads to the notion  \emph{$\mathfrak{F}$-segment} below.

    \begin{definition} Given $T\in\N$, a \emph{$\mathfrak{F}$-segment} of  length  $T$ is a sequence $\cE =\{E_{i}\}_{0\leq i\leq T}$ of compact sets
    such that
   \begin{itemize}
   \item  $f(E_{i})=E_{i+1}$,
   \item there is a family $\{D_i\}_{0\leq i\leq T}\subset \mathfrak{F}$ so that  the compact set  $E_{i}$ is contained $D_i$  and $D_T=E_T$
   \end{itemize}
    We call $E_0$ \emph{the entrance of $\cE$} and $E_T$ \emph{the exit of $\cE$}.
    \end{definition}

    \begin{definition} Given two  $\mathfrak{F}$-segments $ \cE=\{E_{i}\}_{0\leq i\leq T}$ and $\cF=\{F_{j}\}_{0\leq j\leq S}$;
    if the exit of $\cE$ contains the entrance of $\cF$, the \emph{concatenation} of $\cE$ and $\cF$ is a $\mathfrak{F}$-segment
    $\cE\star \cF=\{G_{i}\}_{0\leq i\leq T+S}$ defined as follows:
    \begin{displaymath} G_{i}=\left\{\begin{array}{ll}
    f^{i-T}(F_{0})&\textrm{if $i\leq T$}\\
    F_{i-T}&\textrm {if $i>T$}
    \end{array}\right.
    \end{displaymath}
    \end{definition}

   Next straightforward lemma gives a precise meaning to the simple idea that, if one controls the  averages of $\varphi$ along $\mathfrak{F}$ segments, one also controls the
   averages along the concatenation of these segments. This will allow us to build $\mathfrak{F}$-segment of arbitrarily long length on which we control the averages of $\varphi$.
   \begin{lemma}
    Let $\cE_i=\{E_{i,j}\}_{j\in\{0,\cdots, T_i\}}$, $i\in \{0,\cdots, n\}$, be a family of $\mathfrak{F}$-segments of length $T_i$
    so that the exit of $\cE_i$ contains the entrance of $\cE_{i+1}$ for $i\in\{0,\cdots, n-1\}$.
    Denote by  $T=\sum_{i=0}^n T_i$ and let  $\cF=\{F_j\}_{j\in\{0,\dots, T\}}$ be the $\mathfrak{F}$-segment  defined as the concatenation
    $$\cF=\cE_0\star \cE_1\star\cdots\star\cE_n.$$

    Assume that there are $\alpha<\beta$ so that  for any $i$ and any $x\in E_{i,0}$, one has
    $$\frac 1{T_i}\varphi_{_{T_i}}\in[\alpha,\beta].$$

    Then for every $x\in F_0$ one has
    $$\frac 1T\varphi_{_T}\in[\alpha,\beta].$$

   \end{lemma}

Given $x\in D_0\in\mathfrak{F}$,  a positive integer $t$, a point $s=(s_0,s_1,\cdots)\in\{+,-\}^{\N}$ and $T\in\N^{+}\cup\{+\infty\}$,
    we say that \emph{$x$ follows the $t$-pattern $s$ up to time $T$}, if for any $n\in[0,T)$ which is a multiple of $t$,
    we have that $f^{n+1}(x)\in \cup_{D\in\mathfrak{F}^{s_n}}D$.

Given  $t\in\N^{+}$ and $s\in\{+,-\}^{\N}$, we say that \emph{a $\mathfrak{F}$-segment $\cE=\{E_{i}\}_{0\leq i\leq T}$
     follows $t$-pattern $s$}, if for any $x\in E_0$, the point $x$ follows  $t$-pattern $s$ up to time $T$.

    The key lemma in \cite{BBD1} to find zero center Lyapunov exponent set is the following:

\begin{lemma}\label{kl} \cite[Lemma 2.12]{BBD1}Given a flip-flop family $(\mathfrak{F},\varphi)$, we fix two sequences of positive numbers $\{a_{k}\}$ and $\{b_{k}\}$ which will
converge to zero and satisfy that $b_{k}>a_{k}>b_{k+1}$ for any $k\in\N$.

Then there exists a sequence of integers $1=t_{0}<t_{1}<\cdots$ satisfying that
\begin{itemize}
\item each $t_{i+1}$ is a multiple of $t_i$, for any $i\in\N$;
\item   for every integer $k>1$, every member $D\in\mathfrak{F}$ and every pattern $s\in\{+,-\}^{\N}$,
  there exist two integers $T_{+},T_{-}\in\N$ and two $\mathfrak{F}$-segments $\cE^{+}, \cE^{-}$ of lengths $T_{+}$ and $T_{-}$ respectively
  such that:
  \begin{itemize}
  \item[--] the entrances of $\cE^{+}$ and $\cE^{-}$ are contained in $D$;
  \item[--] $T_{+}$ and $T_{-}$ are multiples of $t_{1}$ and satisfy $t_{k-1}<T_{\pm}\leq t_{k}$;
  \item[--] the segments $\cE^{+}$ and $\cE^{-}$ are $(b_{i},t_{i})$-controlled for $i=1,\cdots,k-1$;
  \item[--] for all $x$ in the entrance of $\cE^{+}$ and all $y$ in the entrance of $\cE^{-}$, we have
   $$a_{k}\leq \frac{1}{T_{+}} \varphi_{_{T_{+}}}(x)\leq b_{k} \textrm{ and } -b_{k}\leq \frac{1}{T_{-}} \varphi_{_{T_{-}}}(y)\leq -a_{k};$$
  \item[--] the segments $\cE^{+}$ and $\cE^{-}$ follow the $t_{1}$-pattern $s$.
  \end{itemize}
  \end{itemize}
\end{lemma}

  According to Lemma \ref{kl}, we can find a sequence of orbit  segments $\{x_{i},T_{i}\}$ whose lengths tend to infinity in  the
  sense of time and those segments are $(b_{j},t_{j})$ controlled for any $j<i$. Then the accumulation  $x$ of $x_{i}$ is controlled
  at any scale and finally any  ergodic measure supported on the $\omega$-limit set of $x$ is non-hyperbolic.
   The last item in the Lemma \ref{kl} guarantees that our system has positive topological entropy on $\omega(x)$.

 \subsection{Blender}

   \emph{Blender} is a powerful tool and shows its power in the study of robust non-hyperbolic phenomena.
   There have been several versions of blenders, see for instance \cite{BD1},\cite{BD2}. Recently, \cite{BBD1}
   gives a simplified definition of blender which is also very easy to understand.
   Here, we will state the new definition of blender and use this new definition.

  Before we state the new definition of blender, let's recall some notations in \cite{BBD1}. We denote by $D^{i}(M)$ the set of $C^1$
  embedded $i$-dimensional compact
  discs in compact Riemannian manifold $M$. We endow $D^{i}(M)$ with $C^1$-topology: for any $D\in D^{i}(M)$, which is the image of
  the embedding  $\psi:\mathbb{D}^{i}\mapsto M$
  where $\mathbb{D}^{i}$ is the $i$-dimensional closed unit disc in $\mathbb{R}^{i}$, the $C^{1}$ neighborhood of $D$ is defined as the set of the images of all the embedding maps
  contained in a $C^1$ neighborhood of $\psi$.  For any $D_1, D_2\in D^{i}(M)$, we define the distance
  $$\rho(D_1,D_2)=\ud_{Haus}(TD_1,TD_2)+\ud_{Haus}(T\partial(D_1), T\partial{D_2}),$$
  where $\ud_{Haus}(\cdot,\cdot)$ denotes the Hausdorff distance on the corresponding Grassmann manifold.
  In \cite[Section 3.1]{BBD1}, it is shown that the distance $\rho(\cdot,\cdot)$
  induces the $C^1$-topology in $D^{i}(M)$.

  Let $f\in\diff^{1}(M)$  and $\mathfrak{D}$ be a subset of $D^{i}(M)$. For any $\epsilon>0$, we denote by $\mathcal{V}_{\epsilon}(\mathfrak{D})$ the $\epsilon$ open neighborhood
  of $\mathfrak{D}$ for the distance $\rho$.

  One says that $\mathfrak{D}$ is a \emph{strictly invariant family with  strength $\epsilon>0$}, if for any $D\in\mathcal{V}_{\epsilon}(\mathfrak{D})$, the image $f(D)$ contains an element of $\mathfrak{D}$.

  \begin{definition}\label{dynamical blender}$(Dynamical\,\, Blender)$ Given  $f\in \diff^{1}(M)$. A
  hyperbolic basic set $\Lambda$ is called a \emph{dynamically defined cu-blender}
  of uu-index $i$, if the followings are satisfied:
  \begin{itemize}
  \item  there is a  dominated splitting of the form $T_{\Lambda}M=E^{s}\oplus E^{c}\oplus E^{uu}$ over $\La$;
  where $\dim(E^{s})=Ind(\Lambda)$, $\dim(E^{c})>0$ and $\dim(E^{uu})=i$.

      \item there exists a neighborhood $U$ of $\Lambda$ such that $\Lambda=\bigcap_{n\in\mathbb{Z}}f^{n}(U)$ and there exists an
      $f$-strictly invariant continuous
      cone field $\mathcal{C}^{uu}$ of index $i$ defined on $\overline{U}$;
          \item  there is a strictly invariant family $\mathfrak{D}\subset D^{i}(M)$ of discs with strength $\epsilon>0$ such that every disc contained in
          $\mathcal{V}_{\epsilon}(\mathfrak{D})$ is tangent to $\mathcal{C}^{uu}$ and contained in $U$.

            \end{itemize}

The set $U$ is called the  domain of $\Lambda$, $\mathcal{C}^{uu}$ is called  strong unstable cone field of $\Lambda$ and $\mathfrak{D}$ is
called  strictly invariant family of discs.
\end{definition}

 We denote   the cu-blender by $(\La, U,\cC^{uu},\mathfrak{D})$.
  We can also define the \emph{cs-blender} which is a \emph{cu-blender} for the reversed dynamics.

  \begin{definition}$(Geometric\,\, Blender)$ Let $f\in\diff^1(M)$. An $f$-invariant compact set $\Lambda$ is
  called a \emph{geometric cu-blender} of uu-index $i$, if the followings are satisfied:
  \begin{itemize}
\item   $\Lambda$ is uniformly hyperbolic with  u-index strictly larger than $i$;
  \item there exist  an open family $\mathfrak{D}\subset D^i(M)$  and a $C^{1}$ neighborhood $\mathcal{U}$ of $f$
  such that for any $g\in\mathcal{U}$ and any $D\in\mathfrak{D}$ we have that $W^{s}(\Lambda_{g})\cap D\neq\emptyset$,
  where $\Lambda_g$ is the continuation of $\Lambda$.
  \end{itemize}
  The open family $\mathfrak{D}$ is called the \emph{superposition region} of $\Lambda$.
  \end{definition}

    The definition of dynamically defined blender is only associated to one diffeomorphism and the geometric blender can tell us the properties of an open set of diffeomorphisms. Actually we have the following result:
    \begin{lemma}\label{robustness of blender} \cite[Lemma 3.14]{BBD1} Let $(\Lambda,U,\mathcal{C}^{uu},\mathfrak{D})$ be a dynamically defined blender with strictly invariant family
    of strength $\epsilon$. Then there exists a $C^1$ neighborhood $\mathcal{U}$ of $f$ such that for any $g\in\mathcal{U}$, $\Lambda_g$
    is a geometric blender with superposition region $\mathcal{V}_{\epsilon/2}(\mathfrak{D})$; furthermore
    $(\Lambda_g,U,\mathcal{C}^{uu},\mathcal{V}_{\epsilon/2}(\mathfrak{D}))$ is a dynamically defined blender for $g$.
    \end{lemma}
\subsection{Flip-flop configuration}
    In this paper, we focus on a co-index one robust cycle, called \emph{flip-flop configuration}, which is  formed by a cu-blender and a hyperbolic periodic orbit of different index such that the unstable manifold of the periodic orbit ``crosses" the superposition region of the cu-blender, and every disc in the strictly invariant family of the cu-blender could ``cross" the stable manifold of the periodic orbit. To be specific:
    \begin{definition}\label{d.flip flop} Let $(\Lambda,U,\mathcal{C}^{uu},\mathfrak{D})$ be a dynamically defined cu-blender of uu-index $i$ and $q$ be a hyperbolic periodic point of u-index $i$. We say that $\Lambda$ and $q$ form a \emph{flip-flop configuration}, if there exist a disc $\Delta^{u}\subset W^{u}(q)$ and a compact submanifold with boundary $\Delta^{s}\subset W^{s}(q)\cap U$ such that:
    \begin{enumerate}
    \item $\Delta^{u}\in\mathfrak{D}$ and $f^{-n}(\Delta^{u})\cap \overline{U}=\emptyset$, for any $n\in\N^+$;
    \item there exists an integer $N$ such that for any $n>N$, $f^{n}(\Delta^{s})\cap\overline{U}=\emptyset$; Moreover, for any $x\in\Delta^s$, if $f^{j}(x)\notin\overline{U}$ for some $j>0$, then  $\{f^{j+k}(x)\}_{k\in\N}\cap\overline{U}=\emptyset$;
        \item for any $y\in\Delta^{s}$, $T_yW^s(q)\cap\mathcal{C}^{uu}=\{0\}$;
        \item there exist a compact set $K\subset\Delta^{s}$ and a  number $\delta>0$ such that for any disc $D\in\mathfrak{D}$, the disc $D$ intersects $K$ in a point whose distance to $\partial{D}$ is no less than $\delta$.
    \end{enumerate}

    \end{definition}
    A set $V$ is called a \emph{neighborhood of flip-flop configuration} if its interior contains the set  $$\mathcal{O}_q\cup\overline{U}\cup\,\bigcup_{i\geq 0}f^{i}(\Delta^s)\cup\,\bigcup_{i\geq 0}f^{-i}(\Delta^u).$$
    It's shown in \cite[Proposition 4.2]{BBD1} that the existence of flip-flop configuration is a robust property.
    \begin{lemma}\cite[Lemma 4.6] {BBD1}Let $f\in\diff^1(M)$. Assume that there exists a blender $(\Lambda,U,\mathcal{C}^{uu},\mathfrak{D})$ of uu-index $i$ forming a flip-flop configuration with a hyperbolic periodic point $q$ of u-index $i$. Let $V$ be a  compact neighborhood of the flip-flop configuration.

  If $V$ is chosen small enough, then the maximal invariant set of $V$ has a partially hyperbolic splitting of the form $E^{cs}\oplus E^{uu}$, where $\dim(E^{uu})=i$. Moreover, there exists a strictly $Df$-invariant cone field $\mathcal{C}_{V}^{uu}$ over $V$ which extends the cone field $\mathcal{C}^{uu}$ defined in $U$,  and any vector in $\mathcal{C}_{V}^{uu}$ is uniformly expanded by $Df$.
    \end{lemma}
    \begin{definition}[Split flip-flop configuration]Consider a flip-flop configuration formed by a blender $(\Lambda, U, \mathcal{C}^{uu}, \mathfrak{D})$ of uu-index $i$ and a hyperbolic periodic point $q$ of u-index $i$. We say that this configuration is \emph{split} if there exists a small compact neighborhood $V$ of this configuration such that the maximal invariant set of $V$ admits a partially hyperbolic splitting of the form $E^{ss}\oplus E^{c}\oplus E^{uu}$,  where $\dim(E^{ss})=Ind(\Lambda)$ and $\dim(E^{uu})=i$.
    \end{definition}
    The following proposition gives the existence of split flip-flop configuration, whose proof is given in \cite[Section 5.2]{BBD1}.
    \begin{proposition}\label{p.existence of flip flop} \cite{BBD1}
    Let $\cU$ be an open set of diffeomorphisms such that for any $f\in\cU$, there exist two hyperbolic periodic points $p_f, q_f$ of  u-index $i_p>i_q$ respectively, depending continuously on $f$ and in the same chain class $C(p_f, f)$. Then there exists an open and dense subset $\cV$ of $\cU$ such that for any $f\in\cV$ and any $i\in(i_q,i_p]$,  there exists a split flop-flop configuration formed by a dynamically defined cu-blender of uu-index $i-1$ and a hyperbolic periodic orbit of u-index $i-1$.
    \end{proposition}

\subsection{Flip flop configuration and flip flop family}

The following proposition shows that the dynamics on  a flip-flop configuration     induces  a flip-flop family.

    \begin{proposition}\label{prop flip}\cite[Proposition 4.9]{BBD1} Consider a diffeomorphism $f$ exhibiting  a
    flip-flop configuration formed by  a dynamical blender $\Lambda$  and  a hyperbolic   periodic orbit $\cO_q$  . Let $V$ be a partially hyperbolic neighborhood of this flip-flop configuration, $\mathcal{C}_V^{uu}$ be the associated strong unstable cone field in $V$ and $K$ be the maximal invariant set of $f$ in $V$. Assume that $\varphi:V\rightarrow \R$ is a continuous function that is positive on $\Lambda$ and is  negative on the periodic orbit $\mathcal{O}_{q}$.

    Then there exist  an integer $N\geq 1$ and a flip-flop family $\mathfrak{F}$ with respect to the dynamics $f^{N}$ and the function
   $$\varphi_{_{N}}:=\sum_{j=0}^{N-1}\varphi\circ f^{j} \textrm{  defined on  $\cap_{j=0}^{N-1}f^{-j}(V)$}.$$
    Moreover, given any $\epsilon>0$, one can choose the flip-flop family $\mathfrak{F}=\mathfrak{F}^{+}\cup\mathfrak{F}^{-}$ such that $\cup\mathfrak{F}^+$ $($resp.$\cup\mathfrak{F}^-$$)$ is contained in the $\epsilon$-neighborhood of $\Lambda$ $($resp.$\mathcal{O}_{q}$.$)$
    \end{proposition}
    \begin{Remark} If $\varphi$ is obtained by extending  $\log\norm{Df|_{E^{c}}}$ on $K$ continuously to $V$, then the points in the $\omega$-limit set of an orbit which are controlled at any scale have a vanishing center Lyapunov exponent.
    \end{Remark}
    \begin{Remark}\label{r.position of disc} According to  \cite[Section 4.4]{BBD1}, one can choose the flip-flop family $(\mathfrak{F},\varphi_{_{N}})$ such that:
    \begin{itemize}
    \item the discs in $\mathfrak{F}$  are tangent to the strong unstable cone field $\cC^{uu}_V$ and have uniform diameter;
    \item the disc $\Delta^{u}\subset W^{u}(\mathcal{O}_q)$ in the definition of flip-flop configuration contains a disc which is an element of $\mathfrak{F}$.
    \item Denote by $W^{s}_{loc}(\mathcal{O}_q)$ the connected component of $W^{s}(\mathcal{O}_q)\cap V$, which contains $\mathcal{O}_q$. For any $D\in\mathfrak{F}$, one of the followings  is satisfied:
        \begin{itemize}\item[--] $f^{2N}(D)\in\mathfrak{D}$;
        \item[--] $f^{2N}(D)$ intersects  $W^s_{loc}(\mathcal{O}_q)$ transversely;
        \item[--]  $D$ intersects $W^s_{loc}(\mathcal{O}_q)$ transversely.

    \end{itemize}
    \end{itemize}
    \end{Remark}

%% file: existence.tex
  \section{Existence of  periodic orbits which are controlled at any scale: proof of Theorem~\ref{thmC}}\label{s.proof of theorem B}
Let  $(\Lambda,U,\mathcal{C}^{uu},\mathfrak{D})$ be a  dynamically defined blender and
 $\mathcal{O}_{q}$ be a hyperbolic periodic orbit.
 We assume that there are $\De^s\subset W^s(\mathcal{O}_{q})$ and $\De^u\subset W^u(\mathcal{O}_{q})$ so that
 $(\Lambda,U,\mathcal{C}^{uu}, \mathfrak{D},\mathcal{O}_{q}, \De^s,\De^u)$   is  a split-flip-flop configuration.

 Fix a  partially hyperbolic neighborhood $V$ of the split flip-flop configuration so that the maximal invariant set $\tilde \La$
 in the closure $\bar V$ admits a partially hyperbolic splitting $T_{\tilde{\La}}M=E^{s}\oplus E^c\oplus E^{u}$ with $\dim(E^c)=1$.

 Thus $$\log \|Df|_{_{E^c}}\|\colon\tilde \La\to \RR$$
 is a continuous function. We denote by
 $$\varphi\colon M\to \RR$$
 a continuous extension of $\log \|Df|_{_{E^c}}\|$.
 We denote $$\|\varphi\|_{_{C^ 0}}=\max\{|\varphi(x)|, x\in M\}.$$
 We denote by $W^{s}_{loc}(\mathcal{O}_q)$ the connected component  of $W^{s}(\mathcal{O}_q)\cap V$  which contains $\mathcal{O}_q$.
 By Proposition \ref{prop flip} and Remark \ref{r.position of disc}, there exists a flip-flop family $(\mathfrak{F},\varphi_{_{N}})$ for  $f^{N}$
 such that
\begin{itemize}\item the disc $\Delta^{u}\subset W^{u}(\mathcal{O}_q)$ contains  a disc which is an element of $\mathfrak{F}$.
    \item for any $D\in\mathfrak{F}$, one of the followings  is satisfied:
        \begin{itemize}\item[--] $f^{2N}(D)\in\mathfrak{D}$;
        \item[--] $f^{2N}(D)$ intersects  $W^s_{loc}(\mathcal{O}_q)$ transversely;
        \item[--]  $D$ intersects $W^s_{loc}(\mathcal{O}_q)$ transversely.
\end{itemize}
\end{itemize}

Let $\lambda <0$ denote the center Lyapunov exponent of the orbit of $q$.

We fix two sequences of positive numbers $\{a_{k}\}_{k\in\N}$ and $\{b_{k}\}_{k\in\N}$ which
converge to zero and satisfy that $b_{k}>a_{k}>b_{k+1}$  and $3b_k< |\lambda |$,  for any $k\in\N$. Moreover, we require that $b_1$ is much smaller than the expanding rate of $Df$ along the bundle $E^u$.

 Note that, $\mathfrak{F}$ is also a flip flop family for the function $\frac 1N \varphi_N$.
 Then for the flip-flop family $(\mathfrak{F},\frac{1}{N}\varphi_{N})$ and two sequences of positive  numbers $\{a_{k}\}$ and $\{b_{k}\}$,
 there exists a sequence of  integers  $1=t_0<t_1<\cdots<t_n<\cdots$ given by    Lemma~\ref{kl}.

 \vskip 2mm

  \begin{proposition}\label{bridge} With the notation above, fix a point $s\in\{+,-\}^\ZZ$, an integer $k\in\N$,  and $\epsilon>0$.
  Then , there exists a hyperbolic periodic orbit
  $\mathcal{O}_{p}\subset V$ homoclinically related to $\mathcal{O}_{q}$,
  such that
  $$-4b_{k}<\lambda^{c}(\mathcal{O}_{p})<-a_{k}.$$

  Moreover,
  there exist  $0<\tau_{1}<\tau_{2}<\pi(p)$ such that:
  \begin{itemize}
     \item the orbit segment $\{p, \tau_{1}\}$ is contained  in the $\epsilon$ neighborhood of the
     negative orbit $Orb^{-}(\Delta^{u},f)$;
     \item $f^{\tau_{1}}(p)$ follows the $t_1$-pattern $s$ up to $\tau_2-\tau_1$;
     \item the point  $f^{\tau_{1}}(p)$ is $(2b_{j}, t_j, \tau_2-\tau_1)$ controlled with respect to $\frac 1N\varphi_N$ for $j=1,\cdots,k$ ;
 \item  the orbit segment $\{f^{\tau_2}(p),\dots ,f^{\pi(p)}(p)\}$ is contained in the $\epsilon$ neighborhood of the positive orbit
 $Orb^{+}(\Delta^{s}\cup W^s_{loc}(\cO_q),f)$.
 \end{itemize}
 \end{proposition}

Note that, up to replacing the diffeomorphism $f$ by $f^N$ and  the map $\varphi$ by $\frac 1N\varphi_N$,
 we just need to prove the proposition
for the case $N=1$.

\proof
Let $\lambda_{k}= -b_{k}<0 $. Corresponding to  the dominated splitting $T_{\tilde\La}M=(E^s\oplus E^c)\oplus E^u$,
Lemma~\ref{shadow} provides two positive numbers $L_{k}$ and $d_{k}$ so that, for any $d\in (0, d_k]$,
any $\lambda_k$-quasi hyperbolic periodic $d$-pseudo orbit is $L_k\cdot d$-shadowed by a periodic orbit.

Now the proof of Proposition~\ref{bridge} consists in building $\lambda_k$-quasi hyperbolic   periodic $d$-pseudo-orbits.
For building these pseudo-orbits, we first use the flip-flop family $\mathfrak{F}$ for building
arbitrarily large $\mathfrak{F}$-segment $\cE$ for which $\varphi$ is $(b_j,t_j)$-controlled for $j\leq k$. Then
we will extend positively and negatively the orbit of a point in $\cE$ in order to get
arbitrarily close to the periodic point $q$, so that  one gets a closed $d$-pseudo orbit.
By requiring that this pseudo orbit
spends enough time (but not too much) close to $q$, we will get a $\lambda_k$-hyperbolic periodic  pseudo-orbit
shadowed by a periodic orbit with the announced center Lyapunov exponent.

Since the flip flop family $(\mathfrak{F},\varphi)$ (with assumption $N=1$) has been chosen so that  $\Delta^{u}$ contains an element $D^u_0$ of $\mathfrak{F}$,
 by applying Lemma~\ref{kl}   to $D^u_0$, there are an integer  $T_1$ and  a  $\mathfrak{F}$-segment $\cE_{1}$ of length $T_1$ so that
  \begin{itemize}
  \item $ t_{k-1}<T_1\leq t_k$;
   \item  we denote by $D^u_1$ the exit of $\cE_1$, then the entrance $f^{-T_1}(D^u_1)$ of $\cE_1$ is contained in $D^{u}_0$,
   \item any point $y\in f^{-T_1}(D^u_1)$ is $(b_j,t_j)$ controlled for $j\in\{1,\dots,k-1\}$;
   \item any point $y\in f^{-T_1}(D^u_1)$  satisfies the inequality
   $$-b_{k}\leq \frac{1}{T_{1}} \varphi_{_{T_{1}}}(y)\leq -a_{k}.$$
   \item   the $\mathfrak{F}$-segment $\cE_1$ follows the $t_{1}$-pattern $s$.
  \end{itemize}

 We build inductively a  sequence  $\{T_i\}$ of integers and a sequence   $\{\cE_i\}_{i\geq 1}$ of $\mathfrak{F}$-segments of length $T_i$.
 $T_i$ and $\cE_i$ are obtained by applying Lemma~\ref{kl} to the exit
 $D^u_{i-1}$ of $\cE_{i-1}$  and have  the following properties:
  \begin{itemize}
  \item $t_{k-1}<T_i\leq t_k$;
   \item we denote by $D^u_i\in\mathfrak{F}$ the exit of $\cE_i$, then  the entrance $f^{-T_i}(D^u_i)$ of $\cE_i$
   is contained in $D^{u}_{i-1}$,
   \item any point $y\in f^{-T_i}(D^u_i)$ is $(b_j,t_j)$ controlled for $j\in\{1,\cdots,k-1\}$;
   \item any point $y\in f^{-T_i}(D^u_i)$  satisfies the inequality $$-b_{k}\leq \frac{1}{T_{i}} \varphi_{_{T_{i}}}(y)\leq -a_{k}.$$
   \item  the  $\mathfrak{F}$-segment $\cE_1\star \cdots \star\cE_i$ follows the $t_{1}$-pattern $s$.
  \end{itemize}

Thus $\cE_1\star \cdots \star\cE_n$, for $n\to\infty$, is the arbitrarily large $\mathfrak{F}$-segment
where $\varphi$ is controlled. Note that $D^u_n$ is the exit of $\cE_2\star \cdots \star\cE_n$.
We will now choose a point in the exit $D^u_n$ whose orbit can be extended
positively and negatively in order to get arbitrarily close to $q$.

Consider a positive number $d<d_{k}$ such that $L_{k}\cdot d$ is smaller than $\epsilon$. The number $d$ needs
to be chosen very small and its precise value will be fixed at the end of the proof. All the constructions below depend  on
the choice of $d$.

According to Remark~\ref{r.position of disc} (and the fact that we assume $N=1$), the exit $f^2(D^u_n)$ intersects
  $\Delta^{s}\cup W^s_{loc}(\mathcal{O}_q)$ transversely.

Since $\Delta^{s}\subset W^{s}(\mathcal{O}_{q})$ is compact, there exists an integer $N_{d}\in\N$ such that
 $f^{N_{d}-2}(\Delta^{s})$ and $f^{N_d-2}(W^s_{loc}(\mathcal{O}_q))$ are contained in $W^s_{d/2}(\mathcal{O}_{q})$.

Hence, there exists a point $y_n$ in $D^u_n$ whose positive orbit remains in $V$ and such that
$$ f^{N_d}(y_n)\in W^s_{d/2}(\mathcal{O}_{q})$$

  We denote
$$x_n=f^{-(T_1+\cdots+T_n)}(y_n).$$
Thus $x_n$ is a point  in the entrance of the $\mathfrak{F}$-segment $\cE_1\star \cdots\star \cE_n$, contained in
$D^u_0\subset \De^u$.

Up to increasing $N_d$ if necessary, we may assume that  :

  $$f^{-\ell}(\Delta^{u})\subset W^s_{d/2}(\mathcal{O}_{q}), \textrm{ for every $\ell\geq N_d$}.$$

Let us denote by $x_{n,\ell}=f^{-\ell}(x_n)$ and $\sigma_{n,\ell}$ the orbit segment
$$\sigma_{n,\ell}=\{x_{n,\ell},\cdots,x_n,\cdots,y_n,\cdots , f^{N_d}(y_n)\}.$$
Note that $\sigma_{n,\ell}$ is contained in the maximal invariant set $\tilde \Lambda$ in $\bar V$.
For any $\ell\geq N_d$,  one has $\ud(x_{n,\ell},f^{N_d}(y_n))< d$ so that the orbit segment
$\sigma_{n,\ell}$ is a closed $d$-pseudo orbit whose period is
$$\pi_{n,\ell}=\ell+\sum_{i=1}^n T_i+ N_d.$$

\begin{lemma}\label{l.hyperbolic} There are $(n,\ell)$ so that:
\begin{itemize}
\item  $$-3 b_k <\frac 1{\pi_{n,\ell}} \sum_{j=0}^{\pi_{n,\ell}-1}\varphi(f^j(x_{n,\ell})) < -b_k$$
 \item  $\sigma_{n,\ell}$ is a $\lambda_k$-quasi hyperbolic $d$-pseudo orbit corresponding to the
 splitting $(E^s\oplus E^c)\oplus E^{u}$ over $\tilde{\Lambda}$.

\end{itemize}

\end{lemma}
\proof We will first prove that there exist  $n$ and $\ell$ arbitrarily large such that the first item is satisfied.

 For simplicity, we assume that $q$ is a fixed point.
  By assumption, we have $3b_{k}<|\lambda|$. We choose  $d$  small enough such
  that for any $z\in B_{d/2}(\cO_q)$, we have that $$|\varphi(z)-\lambda|<\frac 14 a_k.$$

  We denote by $T(n)=\sum_{i=1}^{n}T_i$. For any $\ell\geq N_d$ and any positive integer $n$,
  by the choice of $\sigma_{n,\ell}$, we have that:
  \begin{align*}\sum_{j=0}^{\pi_{n,\ell}-1}\varphi(f^{j}(x_{n,\ell}))
  &=\sum_{j=0}^{\ell-N_d-1}\varphi(f^{j}(x_{n,\ell}))
  +\sum_{j=\ell-N_d}^{\ell-1}\varphi(f^{j}(x_{n,\ell}))
  \\
  &\hspace{5mm}+\sum_{j=\ell}^{\ell+T(n)-1}\varphi(f^{j}(x_{n,\ell}))
  +\sum_{j=\ell+T(n)}^{\pi_{n,\ell}-1}\varphi(f^{j}(x_{n,\ell}))
  \\
  &<(\ell-N_d)(\lambda+\frac{1}{4}a_k)+ N_d\norm{\varphi}_{_{C^0}}
  -T(n)a_k+ N_d\norm{\varphi}_{_{C^0}}
  \\
  &=(\ell-N_d)(\lambda+\frac{1}{4}a_k)
  -T(n)a_k+2N_d\norm{\varphi}_{_{C^0}}
  \end{align*}
  Let $$\Gamma_1=\frac{1}{\pi_{n,\ell}}\big((\ell-N_d)(\lambda+\frac{1}{4}a_k)
  -T(n)a_k+2N_d\norm{\varphi}_{_{C^0}}\big).$$
  On the other hand, we have the following estimate:
   \begin{align*}\sum_{j=0}^{\pi_{n,\ell}-1}\varphi(f^{j}(x_{n,\ell}))
  &=\sum_{j=0}^{\ell-N_d-1}\varphi(f^{j}(x_{n,\ell}))
  +\sum_{j=\ell-N_d}^{\ell-1}\varphi(f^{j}(x_{n,\ell}))
  \\
  &\hspace{5mm}+\sum_{j=\ell}^{\ell+T(n)-1}\varphi(f^{j}(x_{n,\ell}))
  +\sum_{j=\ell+T(n)}^{\pi_{n,\ell}-1}\varphi(f^{j}(x_{n,\ell}))
  \\
  &>(\ell-N_d)(\lambda-\frac{1}{4}a_k)-N_d\norm{\varphi}_{_{C^0}}
  -T(n)b_k-N_d\norm{\varphi}_{_{C^0}}
  \\
  &=(\ell-N_d)(\lambda-\frac{1}{4}a_k)
  -T(n)b_k-2N_d\norm{\varphi}_{_{C^0}}
  \end{align*}
  Let $$\Gamma_2=\frac{1}{\pi_{n,\ell}}\big((\ell-N_d)(\lambda-\frac{1}{4}a_k)
  -T(n)b_k-2N_d\norm{\varphi}_{_{C^0}}\big).$$
  Hence the average of $\varphi$ along $\sigma_{n,\ell}$ belongs to an interval of length:
  $$\Gamma_1-\Gamma_2=\frac{\ell-N_d}{\pi_{n,\ell}}\cdot\frac{a_k}{2}+\frac{T(n)}{\pi_{n,\ell}}(b_k-a_k)+\frac{4N_d}{\pi_{n,\ell}}
  \norm{\varphi}_{_{C^0}}.$$

  There exists $n_0$ such that for any $n>n_0$ and  any $\ell$, we have that
  $$\frac{4N_d}{\pi_{n,\ell}}\norm{\varphi}_{_{C^0}}<\frac{1}{4}a_k$$
  which implies that
  \begin{align*}\Gamma_1-\Gamma_2&=\frac{\ell-N_d}{\pi_{n,\ell}}\cdot\frac{a_k}{2}+\frac{T(n)}{\pi_{n,\ell}}(b_k-a_k)
  +\frac{4N_d}{\pi_{n,\ell}}\norm{\varphi}_{_{C^0}}
  \\&<\frac{a_k}{2}+b_k-a_k+\frac{4N_d}{\pi_{n,\ell}}\norm{\varphi}_{_{C^0}}
  \\
  &<b_k-\frac{1}{4}a_k.
  \end{align*}

 \begin{claim}\label{c.gamma} There are arbitrarily large $\ell$ and $n$  such that:
  $$\Gamma_2=\frac{\ell-N_d}{\pi_{n,\ell}}(\lambda-\frac{1}{4}a_k)
  -\frac{T(n)}{\pi_{n,\ell}}b_k-\frac{2N_d}{\pi_{n,\ell}}\norm{\varphi}_{_{C^0}}
  \in[-3b_k, -2b_k].$$
  \end{claim}
  \proof
  As $N_d$ is constant and $\pi_{n,\ell}=\ell +T(n)+N_d$,  one can check that, for any positive number $\delta>0$,
  there exist   $n_1$ and $\ell_1$ large such that for any $n>n_1$ and $\ell>\ell_1$, we have that
  $$\left|\Gamma_2-\left(\frac{\ell}{\ell+T(n)}(\lambda-\frac{1}{4}a_k)-\frac{T(n)}{\ell+T(n)}b_k\right)\right|<\delta$$

Hence, to  prove the claim, we only need to require that
$$ \delta-3b_k  <\frac{\ell}{\ell+T(n)}(\lambda-\frac{1}{4}a_k)-\frac{T(n)}{\ell+T(n)}b_k<-\delta-2b_k,$$
which is equivalent to
 $$\frac{b_k+\delta}{-\lambda+\frac 14 a_k-2b_k-\delta}<\frac{\ell}{T(n)}<\frac{2b_k-\delta}{-3b_k+\delta-\lambda+\frac 14 a_k}.$$
Remember that $0<-\lambda-3b_k<-\lambda-2b_k$.  Thus, when $\delta$ is small, we have that
   $$\frac{b_k+\delta}{-\lambda+\frac 14 a_k-2b_k-\delta}<\frac{2b_k-\delta}{-3b_k+\delta-\lambda+\frac 14 a_k}.$$
  Hence, we can choose $\ell$ and $n$ arbitrarily large such that
  $$\frac{b_k+\delta}{-\lambda+\frac 14 a_k-2b_k-\delta}<\frac{\ell}{T(n)}<\frac{2b_k-\delta}{-3b_k+\delta-\lambda+\frac 14 a_k}.$$
  \endproof

  Combining   Claim~\ref{c.gamma} with the fact that $\Gamma_1-\Gamma_2<b_k-\frac{1}{4}a_k$, we have that
\begin{equation}\label{e.average}
  \frac 1{\pi_{n,\ell}} \sum_{j=0}^{\pi_{n,\ell}-1}\varphi(f^j(x_{n,\ell}))\in (-3b_k,-b_k-\frac{1}{4}a_k)\subset(-3b_k,-b_k).
\end{equation}
  This ends the proof of the first item of Lemma~\ref{l.hyperbolic} and
   it remains to prove that $\sigma_{n,\ell}$ is a $\lambda_k$-hyperbolic string.

 By the choice of $x_n$, we  have that $Orb(x_n,f)\subset \tilde{\Lambda}$.
 Recall that  we have the partially hyperbolic splitting $T_{\tilde{\Lambda}}M=E^{s}\oplus E^{c}\oplus E^{u}$ and the expanding rate in
 the bundle $E^u$ is  much larger than $-\lambda_k$.

To prove  that $\sigma_{n,\ell}$ is a hyperbolic string (for a good choice of $n$ and $\ell$),
 we only need to show that

 \begin{equation}\label{e.hyperbolic}
 \frac{1}{j}\sum_{m=0}^{j-1}\varphi(f^{j}(x_{n,\ell}))<-b_k, \textrm{ for any $j=1,\cdots,\pi_{n,\ell}$}.
 \end{equation}

Since $n$  can be chosen arbitrarily large and $t_k$ is  a constant, we can require that
        \begin{equation}\label{e.nlarge}\frac{N_d+t_k}{\pi_{n,\ell}}\cdot \norm{\varphi}_{_{C^0}}<\frac 14 a_k.
        \end{equation}

\begin{claim}\label{c.pliss} For any $\ell$  large enough and any $j=1, \cdots,\ell$,
        \begin{equation}\label{e.llarge} \frac{1}{j}\sum_{m=0}^{j-1}\varphi(f^{m}(x_{n,\ell}))<\lambda+\frac{1}{2}a_k<-b_k
        \end{equation}
\end{claim}
\proof Since for any point $z\in B_{d/2}(\cO_q)$, we have that
$$|\varphi(z)-\lambda|<\frac 14 a_k.$$
Recall that for any $j\geq N_d$, $f^{-j}(x_n)\in B_{d/2}(\cO_q)$.
Hence, for any $\ell>N_d$, we have that for any $j\leq \ell-N_d$,
$$\frac{1}{j}\sum_{m=0}^{j-1}\varphi(f^{m}(x_{n,\ell}))<\lambda+\frac{1}{4}a_k.$$
For any $\ell-N_d<j\leq \ell$, we have that
\begin{align*}
\frac{1}{j}\sum_{m=0}^{j-1}\varphi(f^{m}(x_{n,\ell}))&<\frac{\ell-N_d}{j}(\lambda+\frac{1}{4}a_k)
+\frac 1j\sum_{m=\ell-N_d}^{j-1}\varphi(f^{m}(x_{n,\ell}))
\\
&<\frac{\ell-N_d}{\ell}(\lambda+\frac{1}{4}a_k)+\frac{N_d}{\ell-N_d}\norm{\varphi}_{_{C^0}}.
\end{align*}
Hence, when $\ell$ is chosen much larger than $N_d$, we have that
$$\frac{\ell-N_d}{\ell}(\lambda+\frac{1}{4}a_k)+\frac{N_d}{\ell-N_d}\norm{\varphi}_{_{C^0}}<\lambda+\frac{1}{2}a_k.$$
This ends the proof of the Claim~\ref{c.pliss}.
\endproof

 We now choose $n$ and $\ell$ large enough so that Equations~\ref{e.nlarge} and~\ref{e.llarge} hold.
Assume (arguing by contradiction)  that  Equation~\ref{e.hyperbolic} does not hold, then by  Claim~\ref{c.pliss},
there exists an integer $m_0\in(\ell, \pi_{n,\ell})$ such that
$$
        \frac{1}{m_0}\sum_{m=0}^{m_0-1}\varphi(f^{m}(x_{n,\ell}))\geq-b_k.
$$

  Then, combining Equation~\ref{e.average} with Equation~\ref{e.nlarge},  one gets that
  $$m_0\in(\ell,\ell+T(n)].$$

        We denote by $T(i)=\sum_{j=1}^{i}T_j$, for any $i=1,\cdots,n$.
        Then there exists $1\leq i_0\leq n$ such that $$\ell+T(i_0-1)< m_0\leq\ell+T(i_0).$$
 Remember  that $T(i_0)\leq T(i_0-1)+t_k$ and the point $x_n$ is $(b_k,t_k)$-controlled on the time segment $[T(i_0),T(n)]$,
 one gets:

\begin{align*}\sum_{j=0}^{\pi_{n,\ell}-1}\varphi(f^{j}(x_{n,\ell}))&=\sum_{j=0}^{m_0-1}\varphi(f^{j}(x_{n,\ell}))+
\sum_{j=m_0}^{\ell+T(i_0)-1}\varphi(f^{j}(x_{n,\ell}))
        \\
        &\hspace{5mm}+\sum_{j=\ell+T(i_0)}^{\ell+T(n)}\varphi(f^{j}(x_{n,\ell}))+
        \sum_{j=\ell+T(n)+1}^{\pi_{n,\ell}-1}\varphi(f^{j}(x_{n,\ell}))
        \\
        &>-m_0\cdot b_k-t_k\norm{\varphi}_{_{C^0}}-(T(n)-T(i_0))b_k-N_d\norm{\varphi}_{_{C^0}}
        \\
        &=-(m_0+T(n)-T(i_0))b_k-(t_k+N_d)\norm{\varphi}_{_{C^0}}\\
        &>-\pi_{n,\ell}\cdot b_k -(t_k+N_d)\|\varphi\|_{_{C^0}}.
        \end{align*}

        Using Equation~\ref{e.nlarge}, one  gets that
        $$\frac 1\pi_{n,\ell}\sum_{j=0}^{\pi_{n,\ell}-1}\varphi(f^{j}(x_{n,\ell}))>-b_k-\frac 14 a_k.$$
        which is a contradiction to Equation~\ref{e.average}.

        Hence $\sigma_{n,\ell}$  is a $\lambda_{k}$-quasi hyperbolic string corresponding to the splitting $(E^{s}\oplus E^{c})\oplus E^{u}$,
        ending the proof of Lemma~\ref{l.hyperbolic}.
  \endproof

       By Lemma \ref{shadow}, we get a periodic point $p$ of period $\pi_{n,\ell}$ such that for any $i=0,\cdots,\pi_{n,\ell}-1$,
      one has  $$\ud(f^{i}(p),f^{i}(x_{n,\ell}))<L_k\cdot d.$$
       Since $L_k\cdot d<\epsilon$, we have that
       \begin{itemize}
       \item the orbit segment $\{p, \ell-1\}$ is in  $\epsilon$ neighborhood of $Orb^{-}(\Delta^{u},f)$;
       \item  the orbit segment $\{f^{-N_d}(p), N_d\}$ is in $\epsilon$ neighborhood of $Orb^{+}(\Delta^{s},f)$.
        \end{itemize}
        When $d$ is small enough, by the uniformly continuous property of $\varphi$, we have that
        \begin{itemize}
         \item $-4b_{k}<\lambda^{c}(\cO_p)<-a_{k}$;
        \item $p$ is  a $( -a_k , E^s\oplus E^c)$  hyperbolic time in the $(L_{k}+1)\cdot d$ neighborhood of $\mathcal{O}_{q}$;
        \item $f^{\ell}(p)$ follows the $t_1$-pattern $s$ until the time $T(n)$;
         \item  the point  $f^{\ell}(p)$ is $(2b_{j}, t_j, \sum_{i=1}^{n}T_i)$ controlled for $j=1,\ldots,k$;
       \end{itemize}
        By Lemma \ref{lp}, the point $p$ has uniform size of stable manifold of dimension $\dim(E^s\oplus E^c)$.
        By the  the domination of the splitting  $(E^s\oplus E^c)\oplus E^{u}$ and uniform expansion of $E^u$, we have that $\mathcal{O}_p$ is homoclinically related to $\mathcal{O}_{q}$ when $d$ is chosen small.

         Let $\tau_1=\ell$ and $\tau_2=\ell+\sum_{i=1}^{n}T_n$. This gives the proof of Proposition \ref{bridge}.
   \endproof
  With the help of Proposition \ref{bridge}, we now give the proof Theorem \ref{thmC}.
  \proof [Proof of Theorem \ref{thmC}]

  Recall that $\widetilde{\cR\cC}(M)$ is the set of diffeomorphisms with a split flip flop configuration.

  We take $f\in\widetilde{\cR\cC}(M)$.
  Fix a point $s\in\{+,-\}^{\N}$ whose orbit is dense in $\{+,-\}^{\N}$ under the left shift.
   We choose a sequence of positive numbers $\{\epsilon_{k}\}$ which tends to zero. We apply Proposition \ref{bridge} to $(k,\epsilon_k)$,
   then we have a hyperbolic periodic orbit $\gamma_{k}=\mathcal{O}_{p_{k}}$  and two integers  $S_{k}<T_{k}$ such that:
  \begin{itemize} \item $\gamma_k$ is homoclinically related to the orbit of  $q_f$;
  \item $\lambda^c(\cO_{p_k})$ tends to 0;
  \item The orbit segment $\{p_k, S_k\}$ is contained  in the $\epsilon_k$ neighborhood of the
     set  $Orb^{-}(\Delta^{u},f)$;
     \item $f^{S_k}(p_k)$ follows the $t_1$-pattern $s$ up to time $T_k-S_k$;
     \item The point  $f^{S_k}(p_k)$ is $(2b_{j}, t_j, T_k-S_k)$ controlled with respect to  $\frac 1N\varphi_N$ for $j=1,\cdots,k$;
 \item  The orbit segment $\{f^{T_k}(p_k),\cdots ,f^{\pi(p_k)}(p_k)\}$ is contained in the $\epsilon_k$ neighborhood of the set
 $Orb^{+}(\Delta^{s}\cup W^s_{loc}(\cO_{q_f}),f)$.
  \end{itemize}

Up to taking a subsequence of $p_k$, we can assume that $f^{S_k}(p_k)$ converges to a point $x_0$.
Let $$K_f=\cap_{n=1}^{\infty}\overline{\cup_{k=n}^{\infty}\gamma_{k}},$$ then $K_f$ is a compact invariant set and  the orbit of $q_f$ is contained in $K_f$.

         We denote by
          $$\tilde{K}_f=\cap_{n=1}^{\infty}\overline{\cup_{k=n}^{\infty}\{f^{S_{k}}(p_{k}),T_{k}-S_k\}} \textrm{ and } K_f^{\prime}=\cap_{n\in\Z}f^{n}(\tilde{K}_f),$$
           then $\tilde K_f$ is a compact set and $K^{\prime}_f$ is a  compact invariant set.

          \begin{claim}\label{c.proper} $\emptyset\neq K_f^{\prime}\subsetneq K_f$.
          \end{claim}
         \proof Consider the accumulation $x_0$ of the sequence $\{f^{S_{n}}(p_{n})\}_{n=0}^{\infty}$, one has that  $x_0$ is contained in the unstable manifold of $q_f$. Moreover, $Orb^{+}(x_0,f)$ belongs to $\tilde{K}_f$ and by the compactness of $\tilde{K}_f$,  we have $\omega(x_0)\subset\tilde{K}_f$ which implies $\omega(x_0)\subset K^{\prime}_f$. This proves that $K^{\prime}_f$ is non-empty.

           By proving  that $\cO_{q_f}$ is not contained in $K^{\prime}_f$, we show that $K^{\prime}_f\subsetneq K_f$. The proof is by contradiction. Assume that  $\mathcal{O}_{q_f}$ is contained in $K_f^{\prime}$, then there exist   two  sequences  of positive integers $\{m_{i}\}$ and $\{n_i\}$ such that
           \begin{itemize}\item $\lim_{i\rightarrow+\infty}f^{m_i}(p_{n_i})=q_f;$
           \item $m_{i}$ belongs to $[S_{n_{i}}, T_{n_{i}}].$
           \end{itemize}

           If we have a subsequence of $\{|m_{i}-T_{n_i}|\}$ is uniformly bounded from above, then  $q_f$ is controlled at any scale associated to
            $\frac{1}{N}\varphi_{N}$ for the reversed dynamics $f^{-N}$; by Lemma~\ref{l.control at any scale}, $q_f$ has zero center Lyapunov exponent, which  contradicts to the hyperbolicity of $q_f$.

            If $|T_{n_{i}}-m_{i}|$ tends to infinity when $i$ tends to infinity, which implies that $q_f$ is controlled at any scale  associated to $\frac{1}{N}\varphi_{N}$ for $f^{N}$.
            Once again, we get a contradiction.
         \endproof

         \begin{claim}\label{c.zero Lyapunov exponent}Any ergodic measure supported on $K_f^{\prime}$ has zero center Lyapunov exponent.
         \end{claim}
         \proof For any ergodic measure $\mu$ supported on $K_f^{\prime}$, we choose a recurrent point $x\in K_f^{\prime}$ in the basin of $\mu$.
         Arguing as the proof of Claim ~\ref{c.proper},  we get that $x$ is controlled at any scale; by Lemma \ref{l.control at any scale}, we have that $\mu$ has zero center Lyapunov exponent.
         \endproof

         By result of Section 2.5 in \cite{BBD1},  we have that $h_{top}(f|_{K_f^{\prime}})$ is positive and any ergodic measure
         supported on $K_f^{\prime}$ has zero center Lyapunov exponent.

         For any $x\in K_f\backslash (K_f^{\prime}\cup\mathcal{O}_{q_f})$,  there exist two sequences of integers $\{m_i\}$ and $\{n_i\}$
          such that $$\lim_{i\rightarrow+\infty}f^{m_i}(p_{n_i})=x.$$
          Up to taking a subsequence of $m_i$ and $n_i$, we have three possibilities:
         \begin{enumerate}
         \item For each $i$, we have that $m_i\leq S_{n_i}$. Then $x$ belongs to the unstable manifold of $\mathcal{O}_{q_f}$.
          The non-negative number $S_{n_{i}}-m_{i}$ must be uniformly bounded from above.
          Otherwise, $x$ belongs to $\mathcal{O}_{q_f}$ contradicting to our assumption.
           Hence, there exists an integer $N_{x}^{1}$ such that $f^{N_{x}^{1}}(x)$ is an accumulation of $\{f^{S_{n_{i}}}(p_{n_{i}})\}_{i=1}^{\infty}$, which implies that $x$ is controlled at any scale. According to the proof of  Claim \ref{c.proper}, $\omega(x)$ is a subset of $K_f^{\prime}$.
         \item For each $i$, we have that $m_{i}$ belongs to $[T_{n_{i}}, \pi(p_{n_{i}})]$. Then $x$ belongs to $W^{s}(\mathcal{O}_{q})\backslash\mathcal{O}_{q_f}$. Similarly to case one, there exists a number $N_{x}^{2}$ such that $f^{-N_{x}^{2}}(x)$ is accumulation of $\{f^{T_{n_{i}}+1}(p_{n_{i}})\}_{i=1}^{\infty}$, which implies that $x$ is controlled at any scale for map $f^{-1}$. Hence, $\alpha(x)$ is a subset of $K_f^{\prime}$.
         \item For each $i$, we have that $m_{i}$ belongs to $[S_{n_{i}}, T_{n_{i}}]$. Then, we have that either $|m_{i}-S_{n_{i}}|$ or $|T_{n_{i}}-m_{i}|$ is uniformly bounded from above; and we are in the similar situation to case one or case two respectively. Otherwise, according to the definition of $K^{\prime}$, we would have that $\overline{Orb(x,f)}\subset K_f^{\prime}$,  which contradicts to the choice of $x$.
         \end{enumerate}

           This ends the proof of Theorem \ref{thmC}.
         \endproof
\begin{Remark}\label{r.choice of Kf}In the proof of Theorem~\ref{thmC}, one can see that the set $K^{\prime}_f$ contains the $\omega$-limit
set of a point $x_0$ which is controlled at any scale. Actually,  One can choose the $K_f^{\prime}$ to be the $\omega$-limit set of a point in the unstable manifold of $q$, and which
is controlled at any scale.
\end{Remark}
\proof    For the accumulation $x_0$ which is controlled at any scale,  by Proposition \ref{bridge}, there exists a sequence of $\mathfrak{F}$-segments $\{\cD_i\}_{i\in\N}$ such that
         \begin{itemize}
         \item the length of $\cD_i$ tends to infinity;
         \item the entrance of $\cD_i$ is contained in $\Delta^u$;
         \item the entrance of $\cD_i$ tends to the point $x_0$.
         \end{itemize}

     By Remark~\ref{r.position of disc}, all the discs in the flip-flop family has uniform diameter and are  tangent to the strong unstable cone field.
      As a consequence, one has that for any $i\in\N$, there exists a disc $D_i$ in
      $\mathfrak{F}$ such  that the interior of $D_i$ contains $f^i(x_0)$ and $D_i$ is contained in $f^i(\Delta^u)$. Once again, by
     Remark~\ref{r.position of disc}, up to  finite iterates, each disc $D_i $ intersects the local stable manifold of $\cO_q$.
     Now, one can repeat the argument  in   Proposition \ref{bridge} by choosing the quasi hyperbolic pseudo orbit such that
     it spends large proportion
     of time following a long  forward orbit segment  of $x_0$  and spends the rest of time staying close to the local stable and unstable manifolds of $q$. Then using
     the shadowing lemma by \cite{G},  we get a periodic orbit with similar property as the quasi hyperbolic pseudo orbit that we chose.
     As a consequence, up to choosing a subsequence, one gets a sequence of periodic orbits $\{\gamma_n\}_{n\in\N}$ such that
     \begin{itemize}
     \item the center  Lyapunov exponent of $\gamma_n$ tends to zero;
     \item each $\gamma_n$ spends a large proportion of time to follow a long  forward orbit segment of $x_0$ and spends the rest of time
      staying close to  the local stable and unstable manifolds of $q$.
     \end{itemize}

       One can argue as the proof of Theorem~\ref{thmC} above, to show that the compact set $K^{\prime}_f$ is the $\omega$-limit set of $x_0$.

\endproof

    By similar argument, we can have the following result associated to  a sequence of hyperbolic  periodic orbits homoclinically
    related to $\Lambda$ instead of $\cO_{q_f}$.

   \begin{proposition}Assume that $(\Lambda, U, \mathfrak{D},\mathcal{C}^{uu})$ forms a split-flip-flop configuration
   with a hyperbolic periodic orbit $\mathcal{O}_{q}$. Then there exists a sequence of hyperbolic periodic orbits $\{\gamma_{n}\}$ in a neighborhood of the split flip-flop configuration satisfying:
  \begin{itemize} \item $\gamma_{n}$ is homoclinically related to $\Lambda$;
   \item $\lambda^{c}(\gamma_{n})$ tends to $0$;
   \item Consider the set $K=\cap_{n=1}^{\infty}\overline{\cup_{k=n}^{\infty}\gamma_{k}}$, we have that  $\Lambda$ is contained in $K$ and there exists a compact set $K^{\prime}\subset K$ such that  any ergodic measure supported on $K^{\prime}$ has zero center Lyapunov exponent;
       \item $h_{top}(f|_{K^{\prime}})>0$;
       \item  For any $x\in K\backslash (\Lambda\cup K^{\prime})$, we have that either $\omega(x)\subset K^{\prime}$ and $x\in W^{u}(\Lambda)$ or $\alpha(x)\subset K^{\prime}$  and  $x\in W^{s}(\Lambda)$
  \end{itemize}
   \end{proposition}

  \section{Periodic orbits satisfying the \cite{GIKN}  criterion in a flip flop configuration: Proof of Theorem \ref{thmD}}
  Let  $(\Lambda,U,\mathcal{C}^{uu},\mathfrak{D})$ be a  dynamically defined blender and
 $\mathcal{O}_{q}$ be a hyperbolic periodic orbit. Let $\epsilon_0$ be the strength of the strictly invariant family.
 We assume that there are $\De^s\subset W^s(\mathcal{O}_{q})$ and $\De^u\subset W^u(\mathcal{O}_{q})$ so that
 $(\Lambda,U,\mathcal{C}^{uu}, \mathfrak{D},\mathcal{O}_{q}, \De^s,\De^u)$   is  a split-flip-flop configuration.

We fix a  partially hyperbolic neighborhood $V$ of the split flip-flop configuration
so that the maximal invariant set $\tilde \La$ of $f$ in the closure $\bar V$ admits a partially hyperbolic splitting $E^{s}\oplus E^c\oplus E^{u}$ with $\dim(E^c)=1$.
 Let $\varphi\colon M\to \RR$ be the continuous extension of the continuous function $\log \|Df|_{_{E^c}}\|\colon\tilde\La\to \R$.

 Since $\varphi_{|_{\Lambda}}>0$ and  $\Lambda$ is the maximal invariant set of $U$,
 hence  there exist a number $\tau>0$ and an integer $N$ such that for any $x\in \cap_{i=-N}^{N}f^{i}(U)$,
  we have that $$\varphi(x)\geq 2\tau.$$

  \begin{lemma}\label{descend} With the notation above. There exist two constants $\rho\in(0,\frac{1}{\norm{\varphi}_{_{C^0}}})$ and $\zeta\in(0,1)$, such that  for any $\epsilon>0$ and  any hyperbolic periodic orbit $\gamma$  which is contained inside $V$ and is homoclinically related to $\mathcal{O}_{q}$   inside $V$,
    there exists a hyperbolic  periodic orbit $\gamma^{\prime}$ which is homoclinically related to $\gamma$ in $V$
   satisfying:
  \begin{itemize}
  \item $\gamma^{\prime}$ is $(\epsilon, 1-\rho\cdot|\lambda^c(\gamma)|)$ good for $\gamma$;
  \item $\lambda^{c}(\gamma^{\prime})>\zeta\lambda^c(\gamma)$.
  \end{itemize}
  \end{lemma}
  \proof
  We denote by  $\lambda$  the center Lyapunov exponent of $\gamma$,
  then there exists a point $y\in\gamma$ such that $$\frac{1}{k}\sum_{i=0}^{k-1}\log\|Df|_{E^{c}(f^{i}(y))}\|\leq\lambda, \textrm{for $k=1,\cdots,\pi(\gamma)$}.$$

Consider the continuous function $$h_1(t)=\frac{2\norm{\varphi}_{_{C^0}}}{2\norm{\varphi}_{_{C^0}}+|\lambda-t|}\cdot\frac{\lambda+t}{2}+t$$
and
        $$h_2(t)=\frac{2\norm{\varphi}_{_{C^0}}-\tau}{2\norm{\varphi}_{_{C^0}}+|\lambda+t|}\lambda-\frac{3}{2}t,$$
for any $t\geq 0$.

   Since $h_1(0)<\frac{1}{4}\lambda$ and $h_2(0)>\frac{2\norm{\varphi}_{_{C^0}}-\tau}{2\norm{\varphi}_{_{C^0}}}\lambda$,
   there exists $t_0$ such that for any $t\in[0,t_0]$,
    we have the following:
    $$h_1(t)<\frac{\lambda}{4}\textrm{ and } h_2(t)>\frac{2\norm{\varphi}_{_{C^0}}-\tau}{2\norm{\varphi}_{_{C^0}}}\lambda.$$
We take a small positive number $$\delta<\min\{t_0,\,\frac{1}{100}|\lambda|\}.$$

For the number $ \frac{1}{4}\lambda<0$ and the splitting $T_{\tilde\La}M=(E^s\oplus E^c)\oplus E^u$, by Lemma \ref{shadow},
there exist two numbers $L>1$ and $d_{0}>0$ such that for any $d\in(0,d_0]$,
any $\frac{1}{4}\lambda$-quasi hyperbolic $d$-pseudo orbit is $L\cdot d$ shadowed by a real  orbit.
Now we choose a number $d\in(0,d_{0})$ small enough  such that
  \begin{itemize}\item $(L+1)d<\epsilon$;
  \item $|\varphi(z)-\varphi(w)|<\delta$ for any two points $z,w$ satisfying $z\in B_{L\cdot d}(w)$.
   \end{itemize}
   The precise choice of $d$ would be fixed at the end.

   The proof of Lemma \ref{descend} consists in finding a quasi hyperbolic string which starts at a point on the unstable manifold of $y$,
   whose orbit is contained in $\tilde\La$, such that it spends a very long time to follow the periodic orbit $\gamma$.
    Then it spends some  proportion of time in the open set $U$ to gain some expansion in the center direction
    and after that in a small proportion of time it goes into a small neighborhood of $\cO_q$.
     Using  the fact that $\gamma$ and $\cO_q$ are homoclinically related in $V$,
      by the shadowing lemma for hyperbolic set, we can find a hyperbolic string starts from a small neighborhood of $\cO_q$ to $y$.

  Since $\gamma$ is homoclinically related to $\mathcal{O}_q$ inside $V$, by Inclination Lemma,
 there exists an $i$-dimensional compact disc $D^u\in W^{u}(y)\cap V$ and a positive integer $n_1$ such that
 \begin{itemize}
 \item the backward orbit of $f^{n_1}(D^u)$ is contained in $V$;
  \item $f^{n_1}(D^u)$ is $C^1$ $\epsilon_0/2$-close to $\Delta^u$, which implies that $f^{n_1}(D^u)\in\mathcal{V}_{\epsilon_0}(\mathfrak{D})$.
       \end{itemize}

       We denote by $D^u_0=f^{n_1}(D^u)$.
       By the compactness of $D^u$ and of $\Delta^s$,  there exists an integer $n_d$ such that
        $f^{-n_d}(D^u_0)\subset W^u_{d/2}(y)$ and $f^{n_d}(\Delta^s)\subset W^s_{d/2}(\cO_q)$.

       By shadowing lemma for hyperbolic set,  up to increase $n_d$, there exists an $ \lambda/2 $-quasi hyperbolic string
        $\{w,n_{d}\}$ from  $d/2$-neighborhood of $q$ to $d/2$-neighborhood of $y$.

        By the strictly invariant property of $\mathfrak{D}$,   for any integer $r$, we have that
        \begin{itemize} \item $f^{r}(D_0^u)$ contains a uu-disc $D^{u}_r\in\mathfrak{D}$;
        \item $f^{-i}(D^{u}_r)$ is contained in $U$ for any $i=0,\cdots,r$.
        \end{itemize}

       By the definition of flip-flop configuration,   $D^u_r$ intersects $\Delta^s$  in a point $y_r$ transversely, for any positive integer $r$.
       We denote by $x_r=f^{-r}(y_r)$.
   By the choice of $x_r$, one gets that
   \begin{itemize}
   \item the orbit segment $\{x_r, r\}$ is contained in $U$ and $x_r$ belongs to $\tilde\La$;
   \item  for any $n>n_d $ such that $n-n_d$ is a multiple of $\pi(\gamma)$,
        we have  that  $f^{-n}(x_r)\in W^u_{d/2}(y)$.
   \end{itemize}
       For any $r\geq 2N$, where is $N$ is the integer fixed at the beginning of this section, one has that
        \begin{align*}\sum_{i=0}^{r-1}\varphi(f^{i}(x_r))
        &=\sum_{i=0}^{N-1}\varphi(f^i(x_r))+\sum_{i=N}^{r-N-1}\varphi(f^i(x_r))+\sum_{i=r-N}^{r-1}\varphi(f^i(x_r))
        \\
        &> -N\norm{\varphi}_{_{C^0}}+(r-2N)\tau-N\norm{\varphi}_{_{C^0}}
        \\
        &\geq r\cdot\tau-4N\norm{\varphi}_{_{C^0}}
       \end{align*}

       Denote by $x_{r,n}=f^{-n}(x_r)$ and $\sigma_{n,r}$ the orbit  segment
        $$\sigma_{n,r}=\{x_{r,n},\cdots,x_r,\cdots, y_r,\cdots, f^{n_d}(y_r)\}$$
        which is contained in $\tilde\La$. We denote by $\pi_{n,r}=n+r+n_d$.

       \begin{claim}\label{choice of rn}
       There exist two integers $n$ and $r$ arbitrarily large such that
       \begin{itemize}
       \item  $$\frac{n-n_d}{\pi_{n,r}+n_d}>1-\frac{2|\lambda|}{3\norm{\varphi}_{_{C^0}}};$$
       \item $$\frac{2\norm{\varphi}_{_{C^0}}-\tau}{2\norm{\varphi}_{_{C^0}}}\lambda
           <\frac{1}{\pi_{n,r}}\sum_{i=0}^{\pi_{n,r}-1}\varphi(f^{i}(x_{r,n}))<\frac{1}{4}\lambda;$$
       \item $\sigma_{n,r}$ is a $\frac{1}{4}\lambda$-quasi hyperbolic string corresponding to the splitting $(E^s\oplus E^c)\oplus E^{u}$.
       \end{itemize}
       \end{claim}

       \proof By the choice of $\sigma_{n,r}$, we have that
       \begin{align*}\sum_{i=0}^{\pi_{n,r}-1}\varphi(f^{i}(x_{r,n}))&=\sum_{i=0}^{n-n_d-1}\varphi(f^{i}(x_{r,n}))+
       \sum_{i=n-n_d}^{\pi_{n,r}-1}\varphi(f^{i}(x_{r,n}))
       \\
       &<(n-n_d)(\lambda+\delta)+(r+2n_d)\norm{\varphi}_{_{C^0}}.
       \end{align*}

       On the other hand,  we have that
       \begin{align*}\sum_{i=0}^{\pi_{n,r}-1}\varphi(f^{i}(x_{r,n}))&=\sum_{i=0}^{n-n_d-1}\varphi(f^{i}(x_{r,n}))+
       \sum_{i=n-n_d}^{n-1}\varphi(f^{i}(x_{r,n}))
       \\
       &\hspace{5mm}+\sum_{i=n}^{n+r-1}\varphi(f^{i}(x_{r,n}))+\sum_{i=n+r}^{\pi_{n,r}-1}\varphi(f^{i}(x_{r,n}))
       \\
       &>(n-n_d)(\lambda-\delta)-n_d\norm{\varphi}_{_{C^0}}
       \\
       &\hspace{5mm}+ r\cdot\tau-4N\norm{\varphi}_{_{C^0}}-n_d\norm{\varphi}_{_{C^0}}
       \\
       &=(n-n_d)(\lambda-\delta)+r\cdot\tau-(2n_d+4N)\norm{\varphi}_{_{C^0}}.
       \end{align*}
       Hence, there exists $N_0$ such that for any integer $n>N_0$ and any  $r\in\N$, we have that
       $$
       \frac{1}{\pi_{n,r}}\Big((n-n_d)(\lambda+\delta)+(r+2n_d)\norm{\varphi}_{_{C^0}}\Big)
       <\frac{n}{n+r}\lambda
       +\frac{r}{n+r}\norm{\varphi}_{_{C^0}}+\delta$$
       and $$\frac{1}{\pi_{n,r}}\Big((n-n_d)(\lambda-\delta)+r\cdot\tau-(2n_d+4N)\norm{\varphi}_{_{C^0}}\Big)>\frac{n}{n+r}\lambda
       +\frac{r}{n+r}\tau-\delta.$$

       There exist  $n$ and $r$ arbitrarily large such that
       $$\frac{r}{n}\in\Big(\frac{|\lambda+\delta|}{2\norm{\varphi}_{_{C^0}}},
       \frac{|\lambda-\delta|}{2\norm{\varphi}_{_{C^0}}}\Big).$$
       Hence, we have the following estimate:
       \begin{align*}\frac{n}{n+r}\lambda
       +\frac{r}{n+r}\norm{\varphi}_{_{C^0}}+\delta &=\frac{1}{1+\frac rn}\lambda
       +\frac{\frac rn}{1+\frac rn}\norm{\varphi}_{_{C^0}}+\delta
       \\
       &<\frac{1}{1+\frac{|\lambda-\delta|}{2\norm{\varphi}_{_{C^0}}}}\lambda+
       \frac{\frac{|\lambda-\delta|}{2\norm{\varphi}_{_{C^0}}}}{1+\frac{|\lambda-\delta|}{2\norm{\varphi}_{_{C^0}}}}
       \norm{\varphi}_{_{C^0}}+\delta
       \\
       &=\frac{1}{1+\frac{|\lambda-\delta|}{2\norm{\varphi}_{_{C^0}}}}\cdot\frac{\lambda+\delta}{2}+\delta
       \\
       &=\frac{2\norm{\varphi}_{_{C^0}}}{2\norm{\varphi}_{_{C^0}}+|\lambda-\delta|}\cdot\frac{\lambda+\delta}{2}+\delta
       \\
       &=h_1(\delta)
       \end{align*}

       and
       \begin{align*}
       \frac{n}{n+r}\lambda
       +\frac{r}{n+r}\tau-\delta
       &=\frac{1}{1+\frac{r}{n}}\lambda+\frac{\frac{r}{n}}{1+\frac{r}{n}}\tau-\delta
        \\
       &>\frac{1}{1+\frac{|\lambda+\delta|}{2\norm{\varphi}_{_{C^0}}}}\lambda+
       \frac{\frac{|\lambda+\delta|}{2\norm{\varphi}_{_{C^0}}}}{1+\frac{|\lambda+\delta|}{2\norm{\varphi}_{_{C^0}}}}
       \tau-\delta
       \\
       &=\frac{2\norm{\varphi}_{_{C^0}}}{2\norm{\varphi}_{_{C^0}}+|\lambda+\delta|}\lambda
       -\frac{\lambda+\delta}{2\norm{\varphi}_{_{C^0}}+|\lambda+\delta|}\tau-\delta
       \\
       &>\frac{2\norm{\varphi}_{_{C^0}}-\tau}{2\norm{\varphi}_{_{C^0}}+|\lambda+\delta|}\lambda-\frac{3}{2}\delta
       \\
       &=h_2(\delta)
       \end{align*}
       By the choice of $\delta$, we have that
           $$\frac{2\norm{\varphi}_{_{C^0}}-\tau}{2\norm{\varphi}_{_{C^0}}}\lambda
           <\frac{1}{\pi_{n,r}}\sum_{i=0}^{\pi_{n,r}-1}\varphi(f^{i}(x_{r,n}))<\frac{1}{4}\lambda.$$
           This proves the second item of Claim~\ref{choice of rn}.

 Since we have $$\frac{r}{n}\in\big(\frac{|\lambda+\delta|}{2\norm{\varphi}_{_{C^0}}},
       \frac{|\lambda-\delta|}{2\norm{\varphi}_{_{C^0}}}\big),$$ where $n$ and $r$ can be chosen arbitrarily large;
       when $n$ and $r$ are chosen large enough, we have the following
$$\frac{n-n_d }{\pi_{n,r}+n_d}=\frac{1-\frac{n_d }{n}}{1+\frac{r}{n}+\frac{2n_d}{n}}
>\frac{1-\frac{n_d }{n}}{1+\frac{|\lambda-\delta|}{2\norm{\varphi}_{_{C^0}}}+\frac{2n_d}{n}}
>\frac{1}{1+\frac{2|\lambda|}{3\norm{\varphi}_{_{C^0}}}}
>1-\frac{2|\lambda|}{3\norm{\varphi}_{_{C^0}}}.$$

Since $\sigma_{n,r}$ is contained in $\tilde\La$ and $\tilde\La$ admits the partially hyperbolic splitting $(E^s\oplus E^c)\oplus E^{u}$,
to prove that $\sigma_{n,r}$ is a $\frac{1}{4}\lambda $-quasi hyperbolic string, we only need to show that for any integer $j\in[1,\pi_{n,r}]$,
we have the following
$$\frac 1j\sum_{k=0}^{j-1}\varphi(f^{k}(x_{r,n}))\leq \frac{1}{4}\lambda.$$
For any $j\in[1,n+N]$, when $n$ is chosen large enough and $d$ is small enough, we have  the following:
$$\frac{1}{j}\sum_{k=0}^{j-1}\varphi(f^{k}(x_{r,n}))\leq \frac{1}{2}\lambda<\frac{1}{4}\lambda;$$
For any $j\in[n+N,\pi_{n,r}]$, we have that:
\begin{align*}\frac{1}{j}\sum_{k=0}^{j-1}\varphi(f^{k}(x_{r,n}))&=\frac{1}{j}\sum_{k=0}^{n-n_d-1}\varphi(f^{k}(x_{r,n}))
+\frac{1}{j}\sum_{k=n-n_d}^{j-1}\norm{\varphi}_{_{C^0}}
\\
&<\frac{n-n_d}{j}(\lambda+\delta)+\frac{j-n+n_d}{j}\norm{\varphi}_{_{C^0}}.
\end{align*}
Since the last item of the inequality above is increasing when $j$ increases in $j\in[n+N,\pi_{n,r}]$, one has that
$$\frac{1}{j}\sum_{k=0}^{j-1}\varphi(f^{k}(x_{r,n}))\leq \frac{1}{\pi_{n,r}}\big((n-n_d)(\lambda+\delta)+(r+2n_d)\norm{\varphi}_{_{C^0}}\big).$$
By the proof of item two,  one has that
\begin{align*}
& \frac{1}{\pi_{n,r}}\big((n-n_d)(\lambda+\delta)+(r+2n_d)\norm{\varphi}_{_{C^0}}\big)
\\
 &<\frac{n}{n+r}\lambda
       +\frac{r}{n+r}\norm{\varphi}_{_{C^0}}+\delta
       \\
 &<\frac{1}{4}\lambda.
 \end{align*}
This ends the proof of Claim~\ref{choice of rn}.
\endproof

  By the choice of the $ \lambda/2 $-quasi hyperbolic string $\{w,n_d\}$ and the Claim \ref{choice of rn}, we get a $ \frac{\lambda}{4}$-quasi hyperbolic periodic $d$-pseudo orbit $\{\sigma_{n,r},\{w,n_d\}\}$ of period  $\pi_{n,r}+n_d$.
  By Lemma \ref{shadow}, there exists a periodic orbit point $p$ of period $\pi_{n,r}+n_d$ such that
\begin{itemize}
\item For any $j\in[0,\pi_{n,r}]$, we have that
 $$\ud(f^j(p),f^j(x_{r,n}))<L\cdot d;$$
\item For any $j\in[\pi_{n,r}+1, \pi_{n,r}+n_d-1]$, we have that
$$\ud(f^j(p), f^{j-\pi_{n,r}}(w))<L\cdot d.$$
\end{itemize}
Let $\gamma^{\prime}$ be the orbit of periodic point $p$.
Since $(L+1)\cdot d<\epsilon$, we have that $\gamma^{\prime}$ is $(\epsilon, \frac{n-n_d }{\pi_{n,r}+n_d})$ good for $\gamma$.

We take
$$\rho=\frac{2}{3\norm{\varphi}_{_{C^0}}}<\frac{1}{\norm{\varphi}_{_{C^0}}}\textrm{ and }
\zeta=\frac{2\norm{\varphi}_{_{C^0}}-\tau}{2\norm{\varphi}_{_{C^0}}}<1.$$

 By the first item of  Claim \ref{choice of rn}, we have that $\gamma^{\prime}$ is $(\epsilon, 1-\rho\cdot|\lambda|)$ good for $\gamma$.

When $d$ is chosen small enough, by the uniform continuity of $\varphi$ and the third item of Claim \ref{choice of rn},  one gets that

 $$\frac 1j\sum_{k=0}^{j-1}\varphi(f^{k}(p))\leq \frac{1}{5}\lambda,\textrm{ for any integer $j\in[0,\pi_{n,r}+n_d-1]$;}$$
by the uniform continuity of $\varphi$ and the second item of Claim ~\ref{choice of rn}, one gets that $$\lambda^c(\gamma^{\prime})=\frac{1}{\pi_{n,r}+n_d}\sum_{j=0}^{\pi_{n,r}+n_d-1}\varphi(f^j(p))
     >\zeta\lambda.$$

Hence, $p$ is a $( \frac{\lambda}{5} , E^s\oplus E^c)$ hyperbolic time  whose distance to $y$ is less than $(L+1)\cdot d$.
By Lemma \ref{lp},    $p$ has uniform size of stable manifold of dimension $\dim(E^s\oplus E^c)$. When $d$ is small enough, combining
with the fact that $E^{u}$ is uniformly expanding, we have that $\gamma^{\prime}$ is homoclinically related to $\gamma$ in $V$.

This ends the proof of Lemma \ref{descend}.

  \endproof
  \begin{Remark}\label{r.global} If $f$ is globally partially hyperbolic with center dimension one, we can see from the proof of  Lemma \ref{descend} that we can take $V$ as the whole manifold $M$.
  \end{Remark}

Now we can give the proof of Theorem \ref{thmD}.
\proof  We fix a sequence of positive numbers $\{\epsilon_{i}\}$ such  that
 $$\lim_{n\rightarrow +\infty}\sum_{i=0}^{n}\epsilon_{i}<+\infty.$$

Using Lemma \ref{descend}, we will inductively find a sequence of hyperbolic periodic orbits satisfying the condition of Lemma \ref{limit}.

Let $\rho\in(0, \frac{1}{\norm{\varphi}_{_{C^0}}})$ and $\zeta\in(0,1)$ be the two numbers given by Lemma \ref{descend}.

We denote by $$\gamma_{0}=\mathcal{O}_{q} \textrm{ and }\kappa_0=1-\rho\cdot |\lambda^c(\gamma_0)|.$$
  Assume that we already get $\gamma_n$. Then we denote by $$\kappa_n=1-\rho\cdot |\lambda^c(\gamma_n)|.$$
 By applying  Lemma \ref{descend} to $\gamma_n$ and $\epsilon_n$, we get a hyperbolic periodic orbit $\gamma_{n+1}$ such that
  \begin{itemize}\item $\gamma_{n+1}$ is homoclinically related to $\gamma_n$  in $V$ ;
  \item $\gamma_{n+1}$ is $(\epsilon_n, \kappa_n)$ good for $\gamma_n$;
  \item $|\lambda^c(\gamma_{n+1})|<\zeta \cdot |\lambda^c(\gamma_n)|$.
  \end{itemize}
 For any $n$, we have that
 $$|\lambda^c(\gamma_n)|\leq \zeta^n\cdot |\lambda^c(\gamma_0)|.$$
 Hence, the center Lyapunov exponent of $\gamma_n$ exponentially tends to zero when $n$ tends to infinity,
  which implies  $$\lim_{n\rightarrow +\infty}\prod_{i=0}^{n}\kappa_i\in(0,1].$$

 By Lemma \ref{limit},
   the sequence $\{\delta_{\gamma_{n}}\}$ converges to a non-hyperbolic ergodic measure $\mu$ whose support is $$\cap_{n=1}^{\infty}\overline{\cup_{k=n}^{\infty}\gamma_{k}}.$$
\endproof 

%% file: proofs.tex
\section{Existence of non-hyperbolic ergodic measures with full support for robustly transitive  diffeomorphisms:
 Proof of Theorem \ref{thm.non-hyperbolic ergodic measure} }

Recall that $\cT(M)$ is the set of robustly transitive  partially hyperbolic (but non-hyperbolic) diffeomorphisms,
whose center can  be split into one dimensional subbundles which form dominated splittings.
We denote by $d=\dim(M)$ and  $i_0=\dim(E^{s})$.

Up to changing a metric (due to \cite{Go}), we can assume that there exists $\lambda_0<0$ such that
\begin{itemize}\item $\log\norm{Df|_{E^s(x)}}<\lambda_0$ and $\log\norm{Df^{-1}|_{E^{u}(x)}}<\lambda_0$, for any $x\in M$;
\item  For any $x\in M$, we have  that
$$\log\norm{Df|_{E^s(x)}}-\log\norm{Df|_{E^c_1(x)}}<2\lambda_0,$$
$$\log\norm{Df^{-1}|_{E^{u}(f(x))}}+\log\norm{Df|_{E^c_k(x)}}<2\lambda_0,$$
and  for any $i=1,\cdots,k-1$,
$$\log\norm{Df|_{E^c_{i}(x)}}-\log \norm{Df|_{E^c_{i+1}(x)}}<2\lambda_0.$$
\end{itemize}

Let $p_f$ be a $f$-hyperbolic periodic point. We say that the homoclinic class of $p_f$ is \emph{robustly being the whole manifold}, if there exists a $C^1$ small neighborhood $\cU_f$ of $f$ such that for any $g\in\cU_f$, we have that
\begin{itemize}
\item  the continuation $p_g$ of $p_f$ is well defined;
\item  the homoclinic class of $p_g$ is the whole manifold.
\end{itemize}
\subsection{Existence of homoclinic classes robustly being the whole manifold}



By \cite{BC},  for   $C^1$ generic diffeomorphisms in $\cT(M)$ and   any $j=0,\cdots, k$,  the set of periodic orbits of index $i_0+j$ is  dense on $M$ and periodic orbits of the same index are homoclinically related; As a consequence, we have that $M$ is a homoclinic class. Recently, \cite{ACS} proves that one can replace the generic assumption by open and dense assumption to show that $M$ is a homoclinic class of periodic orbits of index $i_0$ and $i_0+k$ in a robust way. Combining with \cite[Theorem E]{BDPR}, we have the following:

\begin{proposition}\label{p.homoclinic class}There exists an open and dense subset $\cT_h(M)$ of $\cT(M)$ such that for any $f\in\cT_h(M)$,  there exist $k+1$ hyperbolic periodic points $p_1,\cdots, p_{k+1}$ whose homoclinic classes are robustly being the whole manifold.
\end{proposition}

\subsection{Periodic orbits satisfying \cite{GIKN}    criterion: Proof of Theorem \ref{thm.non-hyperbolic ergodic measure}}
To prove Theorem \ref{thm.non-hyperbolic ergodic measure}, we need the following proposition.

\begin{proposition}\label{p.non-hyperbolic ergodic measure}Let $f\in\cT(M)$ and $p$ be a hyperbolic periodic point  of index $i_0+i$ for some integer $i\in(0, k]$. Assume that the homoclinic class of $p$  is the whole manifold. Assume, in addition, that there exists a cu-blender $(\La^{u},U,\cC^{uu},\mathfrak{D})$ of uu-index $d-i_0-i$  such that  $\mathcal{O}_p$ and $\La^{u}$ form a split flip-flop configuration.

Then there exist $\rho>0$ and $\zeta\in(0,1)$ such that for any $\epsilon>0$ and    any hyperbolic periodic orbit $p_0$ homoclinically related to $p$ satisfying that $\lambda^c_i(p_0)>\lambda_0$, where $\lambda_0$ is the number we fixed at the beginning of this section,
there exists a hyperbolic periodic point $\mathcal{O}_{p_1}$ such that
\begin{itemize} \item $\mathcal{O}_{p_1}$ is homoclinically related to $\cO_{p_0}$ and is $\epsilon$ dense in $M$;
\item The orbit of $p_1$ is $(\epsilon, 1-\rho\cdot|\lambda^c_{i}(p_0)|)$ good for $p_0$;
\item $\lambda^c_i(p_1)>\zeta\cdot\lambda^c_{i}(p_0)$.
\end{itemize}
\end{proposition}
\proof  Since $Df$ is uniformly expanding in the bundle  $E^c_i|_{\La^u}$,  there exist $\tau_1>\tau_2>1$ and an integer $N$ such that
 for any $x\in \cap_{i=-N}^{N} f^{i}(U)$,
 we have that
$$\tau_2<\norm{Df|_{E^c_{i}(x)}}<\tau_1.$$
For simplicity, we will take $N=1$.

 We denote by $\lambda$ the Lyapunov exponent of  $p_0$ along $E^c_i$. 

We take $\delta\in(0, \frac{-\lambda}{4})$, whose precise value would be fixed at  the end.
 By the uniform continuity of the functions $\log\norm{Df|_{E^{c}_i}}$ and $\log\norm{Df|_{E^{c}_{i+1}}}$,
there exists $\eta>0$ such that for any two points $y, w\in M$ satisfying that $\ud(y,w)<\eta$, we have that
$$|\log\norm{Df|_{E^{c}_i(y)}}-\log\norm{Df|_{E^{c}_i(w)}}|<\frac{\delta}{2}$$
  and $$|\log\norm{Df|_{E^{c}_{i+1}(y)}}-\log\norm{Df|_{E^{c}_{i+1}(w)}}|<\frac{\delta}{2}.$$

Since  $H(p_0,f)=M$, for any $\epsilon\in(0,\eta)$ and any $\kappa\in (0,1)$,  there exists a hyperbolic periodic point $p^{\prime}$ homoclinically related to $p_0$ (therefore homoclinically related to $p$) such that
\begin{itemize}
\item the orbit of $p^{\prime}$ is $\frac{\epsilon}{2}$ dense in $M$;
 \item the orbit of $p^{\prime}$ is $(\frac{\epsilon}{2},\kappa)$ good for the orbit of $p_0$.
 \end{itemize}

Since  $E^c_i$ is one dimension,  we have that
$$\frac{1}{\pi(p_0)} \sum_{j=0}^{\pi(p_0)-1}\log\norm{Df|_{E^c_{i}(f^j(p_0))}}=\lambda.$$
By assumption  that $\lambda\in(\lambda_0, 0)$ and  the domination
 $$\frac{1}{\pi(p_0)}\cdot \sum_{j=0}^{\pi(p_0)-1}\log\norm{Df|_{E^c_{i}(f^j(p_0))}}-\frac{1}{\pi(p_0)} \cdot \sum_{j=0}^{\pi(p_0)-1}\log\norm{Df|_{E^c_{i+1}(f^j(p_0))}}\leq 2\cdot\lambda_0,$$
we have that  $$\frac {1}{\pi(p_0)} \sum_{j=0}^{\pi(p_0)-1}\log\norm{Df|_{E^c_{i+1}(f^j(p_0))}}>-\lambda_0.$$
Hence, by Lemma \ref{l.bi-hyperbolic time}, there exists a $ \lambda+\delta/2 $ bi-hyperbolic time on the orbit $p_0$.
 For notational convenience,  we still denote the $ \lambda+\delta/2 $ bi-hyperbolic time as $p_0$.

By the uniform continuity of $\log\norm{Df|_{E^c_i}}$ and of $\log\norm{Df|_{E^c_{i+1}}}$, when  $\epsilon$ is taken small and  $\kappa$ is close to $1$ enough,  the orbit of $p^{\prime}$ has a   $\lambda+\delta$ bi-hyperbolic time in the $\frac \epsilon 2$ neighborhood of $p_0$. For simplicity, we denote the bi-hyperbolic time as $p^{\prime}$.

By Lemma \ref{shadow}, there exist two numbers $L$ and $d_0$ corresponding to the number  $ \frac{\lambda}{4} $ and to  the splitting
$$TM=(E^s\oplus E_1^c\oplus\cdots\oplus E^c_i)\oplus (E^c_{i+1}\oplus\cdots  \oplus E^{u}).$$

We take $d<\epsilon/2$ small enough such that $L\cdot d$ is much smaller than $\epsilon/2$.
 The precise value of $d$ would be fixed at last.
By the definition of flip-flop configuration and the homoclinic relation between $\cO_{p^{\prime}}$ and $\cO_p$, there exists  an integer $N_d$ such that
 \begin{itemize}\item $f^{N_d}(W^{u}_{d/2}(p^{\prime}))$ contains a uu-disc $D^{u}_0\in\mathfrak{D}$;
 \item For any disc $D\in\mathfrak{D}$, we have that $f^{N_d}(D)$ intersects  $W^{s}_{d/2}(p^{\prime})$ transversely.
 \end{itemize}
 For any positive integer $m$, we have that
$f^{m}(D^{u}_0)$ contains a uu-disc  $D^{u}_m\in\mathfrak{D}$. By the choice of $N_d$, we have that $f^{N_d}(D^{u}_m)$ intersects $W^{s}_{d/2}(p^{\prime})$ in a point $x$.

Then for any integers $n, m\textrm{ and } k$, we consider the following orbit segment:
$$\sigma_{n,m,k}=\{f^{-n\pi(p^{\prime})-m-2N_d}(x), \cdots,
x, \cdots ,f^{k\pi(p^{\prime})}(x)\}.$$
We denote by $$x_{n,m}=f^{-n\pi(p^{\prime})-m-2N_d}(x)\in W^u_{d/2}(p^{\prime}) \textrm{ and } \pi_{n,m,k}=(n+k)\pi(p^{\prime})+m+2N_d.$$
Notice that $\ud(x_{n,m}, f^{k\pi(p^{\prime})}(x))<d$.
\begin{claim}\label{c.choice of nmk} There exist  $n, m, k$ which can be chosen arbitrarily large,  such that when $\delta$ is chosen small enough, we have
\begin{itemize}
\item $$\frac{n\pi(p^{\prime})}{\pi_{n,m,k}}>1-\frac{2|\lambda|}{3\log\tau_1}$$
\item  $$\frac{1}{\pi_{n,m,k}}\sum_{j=0}^{\pi_{n,m,k}-1}\log\norm{Df|_{E^c_{i}(f^j(x_{n,m}))}}\in\big(\frac{2\log\tau_1-\log\tau_2}{2\log\tau_1}\cdot \lambda,\frac{\lambda}{4}\big).$$
\item $\sigma_{n,m,k}$ is a $ \frac{\lambda}{4} $ quasi hyperbolic string corresponding to the splitting
$$TM=(E^s\oplus E_1^c\oplus\cdots\oplus E^c_i)\oplus (E^c_{i+1}\oplus\cdots  \oplus E^{u}).$$
\end{itemize}
\end{claim}
The proof of the Claim~\ref{c.choice of nmk}  is similar to Claim~\ref{choice of rn} in the proof of Lemma \ref{descend}. The   difference  is the last item, since the bundle $E^c_{i+1}$ is not uniformly expanding. However, $E^c_{i+1}$ is uniformly expanding when it is restricted to the neighborhood of the blender.  Here, we explain a little bit about   the uniform contraction of $E^c_{i+1}$ by $Df^{-1}$ from $f^{k\pi(p^{\prime})}(x)$ to
$x_{n, m, k}$.
From the proof of Claim~\ref{choice of rn},  when we choose $\delta$ small enough,  for any integer $k$, we can  choose $n$ and  $m$ arbitrarily large  such that $\frac{k}{n}$ are small enough,  the first item and the second item of the claim are satisfied, and  $x_{n,m,k}$ is $( \frac{\lambda}{4}, E^c_i)$ hyperbolic time  until the point $f^{k\pi(p^{\prime})}(x)$. To make sure $\sigma_{n,m,k}$  is a $ \frac{\lambda}{4} $ quasi hyperbolic string, we only need to show that  $f^{k\pi(p^{\prime})}(x)$ is the $(-\frac{\lambda}{4}, E^{c}_{i+1})$  hyperbolic time until the point $x_{n, m, k}$. To guarantee this, we only need to ask that $k$ is much larger than $2N_d+\pi(p^{\prime})$,  but still much smaller than $n$, and we can do that  is because we have the following fact:
\begin{itemize}
\item$$\lambda_0<\lambda<0$$
\item $$\tau_2<\norm{Df|_{E^c_{i}(x)}}<\tau_1, \textrm{ for any $x\in f(U)\cap U\cap f^{-1}(U)$.} $$
\item
$$\log\norm{Df|_{E^c_{i}(z)}}-\log \norm{Df|_{E^c_{i+1}(z)}}<2\lambda_0, \textrm{ for any $z\in M$.}$$
\end{itemize}

We take $\rho=\frac{2}{3\log\tau_1}$ and $\zeta=\frac{2\tau_1-\log\tau_2}{2\log\tau_1}$.

By Lemma \ref{shadow}, there exists a periodic point $p_1$ such that
$$\ud\big(f^{j}(x_{n,m,k}),f^j(p_1)\big)<L\cdot d<\epsilon/2, \textrm{ for any $j=0,\cdots, \pi_{n,m,k}-1$}.$$

By the choice of $\sigma_{n,m,k}$ and Claim \ref{c.choice of nmk}, when $d$ is chosen small,  we have that
\begin{itemize}
\item $\lambda^c_i(p_1)\in(\zeta\cdot\lambda, \frac{\lambda}{4})$;
\item the orbit of  $p_1$ is $(\epsilon, 1-\rho\cdot|\lambda|)$ good for $p_0$;
\item the orbit of $p_1$ is $\epsilon$ dense in $M$.
\end{itemize}

Once again when $d$ is chosen small enough, by Lemma \ref{lp} and  the uniform continuity of the functions
  $\log\norm{Df|_{E^{c}_i}}$ and $\log\norm{Df|_{E^{c}_{i+1}}}$, we have that $p_1$ has the uniform size of stable and unstable manifolds
  which implies that $p_1$ is homoclinically related to $p^{\prime}$.

  This ends the proof Proposition \ref{p.non-hyperbolic ergodic measure}. \endproof

\proof[Proof of Theorem \ref{thm.non-hyperbolic ergodic measure}]  For any $j\in[1, k]$, by Proposition \ref{p.homoclinic class} and Proposition \ref{p.existence of flip flop}, there exists an open and dense subset $\tilde{\cT}_j (M)$ of  $\cT(M)$, such that  for any $f\in\tilde{\cT}_j(M)$,
\begin{itemize}\item there exists a split flip-flop configuration formed by a dynamically defined cu-blender $\La^u_j$ and a hyperbolic periodic orbit $p_j$ of s-index $i_0+j$;
\item there exists a hyperbolic periodic point $q_j$ of s-index $i_0+j$ whose homoclinic class  is robustly being the whole manifold.
\end{itemize}

By connecting lemma and robust transitivity, we can do arbitrarily $C^1$ small perturbation to make $q_j$ and $p_j$ be homoclinically related. As a consequence, there exists an open and dense subset  $\cT_j(M)$ of $\tilde{\cT}_j(M)$, such that for any $f\in\cT_j(M)$,  there exists a split flip-flop configuration formed by a dynamically defined cu-blender $\La^u_j$ and a hyperbolic periodic orbit $p_j$ of index $i_0+j$; moreover the homoclinic class of $p_j$ is robustly being the whole manifold.

 By Lemma \ref{descend},   there always exists a hyperbolic periodic point $q_j^{\prime}$ homoclinically related to $p_j$ whose Lyapunov exponent along the bundle $E^c_{j}$ is   much larger than $\lambda_0$.

Let $\rho>0$ and $\zeta\in(0,1)$ be the two constants given by Proposition \ref{p.non-hyperbolic ergodic measure}. We fix a sequence of positive numbers $\{\epsilon_n\}_{n\geq 1}$   such that $\sum_n\epsilon_n$ converges.

We will inductively find a sequence of hyperbolic periodic orbits satisfying the condition in Lemma \ref{limit}.

  Denote by  $\gamma_0=\cO_{q_j^{\prime}}$.
  Assume that we already have a periodic orbit $\gamma_n$ such that
  \begin{itemize} \item $\gamma_n$ is homoclinically related to $\gamma_0$ and is  $\epsilon_{n}$ dense in $M$;
  \item $\lambda^c_j(\gamma_n)>\zeta\cdot\lambda^c_j(\gamma_{n-1})$;
  \item $\gamma_n$ is $(\epsilon_n, 1-\rho|\lambda^c_j(\gamma_{n-1})|)$ good for $\gamma_{n-1}$.
  \end{itemize}
By Proposition \ref{p.non-hyperbolic ergodic measure}, there exists a hyperbolic periodic orbit $\gamma_{n+1}$ such that
 \begin{itemize} \item $\gamma_{n+1}$ is homoclinically related to $\gamma_0$ and  $\epsilon_{n+1}$ dense in $M$;
  \item $\lambda^c_j(\gamma_{n+1})>\zeta\cdot\lambda^c_j(\gamma_{n})$;
  \item $\gamma_{n+1}$ is $(\epsilon_n, 1-\rho|\lambda^c_j(\gamma_{n})|)$ good for $\gamma_{n}$.
  \end{itemize}

 One can see that $\lambda^c_j(\gamma_n)$ converges to $ 0$ exponentially which implies that the product $\prod_n(1-\rho|\lambda^c_j(\gamma_{n-1})|)$ converges to a positive number.
By Lemma \ref{limit} and continuity of the function $\log\norm{Df|_{_{E_j^c}}}$,  the Dirac measure $\delta_{\gamma_n}$ converges to a non-hyperbolic ergodic measure $\nu_j$ whose support is given by $\cap_{n=1}^{\infty}\overline{\cup_{k=n}^{\infty}\gamma_k}$.  Since $\epsilon_n$ tends to zero, we have that $supp\, \nu_j=M$.

We take the intersection $\tilde{\cT}(M)=\cap_{j=1}^{k}\cT_j(M)$, which is an open and dense subset of $\cT(M)$.
This ends the proof of Theorem \ref{thm.non-hyperbolic ergodic measure}.
\endproof 

%% file: biblio.tex
\bibliographystyle{plain}